\documentclass[preprint,12pt,authoryear]{elsarticle}
\usepackage[utf8]{inputenc}
\usepackage[noautomatic]{imakeidx}
\usepackage{hyperref}%
\usepackage{hyperref}
\hypersetup{colorlinks,linkcolor=blue}
\makeindex

\usepackage{url}
\usepackage{lipsum}
\usepackage{algorithmic}
\usepackage{amsopn}
\DeclareMathOperator{\diag}{diag}
\usepackage[T1]{fontenc}
\usepackage{lmodern}
\usepackage{subfig,color}
\usepackage[english]{babel}
\usepackage{multirow}
\usepackage{amsmath,amssymb,psfrag}
\usepackage{graphicx}
\usepackage{listings}
\usepackage{paralist}
\usepackage{appendix}
\usepackage{amsfonts}
\usepackage{subfig}
\usepackage{comment} 
\usepackage{enumerate}
\usepackage{cancel}
\usepackage{footnote}
\usepackage{epstopdf}
\usepackage{tikz,tikz-cd}
\usepackage{listings}
\usepackage{pgfplots,pgfplotstable}
\graphicspath{{./Figures/}}

\ifpdf
  \DeclareGraphicsExtensions{.eps,.pdf,.png,.jpg}
\else
  \DeclareGraphicsExtensions{.eps}
\fi

\usetikzlibrary{pgfplots.groupplots}
\usetikzlibrary{positioning}
\usetikzlibrary[shapes,arrows,trees]
\usetikzlibrary{matrix,decorations.pathmorphing}

\renewcommand{\d}{\mathrm{d}}

\newcommand{\D}{\mathrm{D}}

\newcommand{\dt}{\Delta t}

\newcommand{\Lagrangian}{{\mathcal{L}}}

\newcommand{\PETSc}{{\textnormal{PETSc}}}
\newcommand{\TAO}{{\textnormal{TAO}}}
\newcommand{\TS}{{\textnormal{TS}}}

\usepackage[hyphens,spaces,obeyspaces]{trl}
\usepackage{amsmath}

\RequirePackage[normalem]{ulem} 
\RequirePackage{color}
\definecolor{RED}{rgb}{1,0,0}\definecolor{BLUE}{rgb}{0,0,1}   

\newcommand{\SErr}{{\rm wlte}}
\newcommand{\Tol}{{\mathrm{Tol}}}
\newcommand{\TolA}{{{\rm Tol}_{\rm A}}}
\newcommand{\TolR}{{{\rm Tol}_{\rm R}}}
\pgfplotsset{compat=1.13}

\journal{Journal of Computational Science}

\begin{document}

\printindex

\numberwithin{equation}{section}
\begin{frontmatter}

\title{A scalable matrix-free spectral element approach for unsteady PDE constrained
optimization using PETSc/TAO}

\author{Oana Marin, Emil Constantinescu, Barry Smith}
\ead{oanam@mcs.anl.gov,emconsta@mcs.anl.gov,bsmith@mcs.anl.gov}

\address{Mathematics and Computer Science Division \\ Argonne National Laboratory \\ Lemont, IL 60439}

\begin{abstract}
  We provide a new approach for the efficient matrix-free application
  of the transpose of the Jacobian for the spectral element method for
  the adjoint based solution of partial differential equation (PDE)
  constrained optimization.  This results in optimizations of nonlinear PDEs
  using explicit integrators where the integration of the adjoint problem
  is not more expensive than the forward simulation.  Solving
  PDE constrained optimization problems  entails combining
  expertise from multiple areas, including simulation, computation of
  derivatives, and optimization. The Portable, Extensible Toolkit for
  Scientific computation (\PETSc{}) together with its companion
  package, the Toolkit for Advanced Optimization (\TAO{}), is an
  integrated numerical software library that contains an
  algorithmic/software stack for solving linear systems, nonlinear
  systems, ordinary differential equations, differential algebraic
  equations, and large-scale optimization problems and, as such, is
  an ideal tool for performing PDE-constrained optimization. This
  paper describes an efficient approach in which the software stack
  provided by \PETSc{}/\TAO{} can be used for large-scale nonlinear
  time-dependent problems.  While time integration can involve a
  range of high-order methods, both implicit and explicit. The
  PDE-constrained optimization algorithm used is gradient-based and
  seamlessly integrated with the simulation of the physical problem.
\end{abstract}


\end{frontmatter}

\section{Introduction}
Fitting numerical or
experimental observations to determine parameters, identifying boundary conditions that satisfy
certain observations, optimizing an objective function of a simulation
solution, accelerating simulations through their long transients (the
spin-up problem \citep{isik2013spin}), and many more computations
fall within the field of partial differential equation (PDE)-constrained optimization and inverse problems.  Despite their
widespread utility, the solution of such problems is plagued by  bottlenecks ranging from mathematical issues to high computational
costs, high input-output costs, and especially software complexity. Several codes, such as JuMP \citep{dunning2017jump}, and 
Python libraries, such as the dolfin-adjoint project, \citep{farrell2013automated}, address these issues; however, they exhibit
limitations for large-scale problems.

A unique contribution of this paper is the efficient, scalable application of
the transpose of the Jacobian matrix-vector product for spectral
elements at a cost not much more than the forward application of the
nonlinear operator using tensor contractions.  The matrix-free
approach avoids the cost and storage required by forming the Jacobian
explicitly. Our implementation of the optimization process leverages the extensive range of time integrators
available in the Portable, Extensible Toolkit for Scientific
in (\PETSc{}) \citep{petsc-user-ref} and the  optimization
algorithms readily available in the Toolkit for
Advanced Optimization (\TAO{}) \citep{tao-user-ref}. In particular, we utilize the adjoint integrators introduced in \citep{Zhang_2019}.


Time-dependent PDE-constrained optimization problems can be posed in two ways: either
by fully discretizing, in both space and time, the entire Karush-Kuhn-Tucker (KKT) system (that
is, solving for all unknowns at all time steps simultaneously, sometimes called the {\em
  all-at-once} approach)  \citep{Haber_2007} or by decoupling the direct problem and
its adjoint \citep{giles1997adjoint,ou_2011,gunzburger2002perspectives,sandu2006properties,gander2014constrained} 
and solving the optimization problem via a forward/backward
time-stepping loop, with appropriate initial and boundary conditions
depending on the objective function. This work focuses on the latter, in particular on discrete
adjoint approaches, targeting a minimal computational cost implementation of the adjoint equation.
Optimization algorithms that involve time integration are computationally expensive, for example, for an optimization that requires $m$ iterations,
one needs to perform $2m$ integrations each of $N$ time steps, where the backward integration can be
more expensive than the forward integration. Employing
all possible computational accelerations is crucial since any speedup accumulates significantly with the number of time steps and optimization iterations.

The core strategy presented in this study is to solve the inverse (or PDE-constrained optimization)
problem in a purely nonlinear manner. Considerable work has previously
relied on {\em linearized adjoints}; that is, the model is first
linearized, and the adjoints are computed for that linearized problem
\citep{saglietti2016adjoint}. Here we showcase the solution of
the nonlinear time-dependent viscous Burgers equation utilizing
\PETSc{} with \TAO{}.

In previous work \citep{schanen_2016}, we addressed the often-overlooked data intensity and
computational intensity aspects of PDE-constrained optimization. The data intensity for nonlinear problems stems
from the nature of the backward-in-time component, which depends at
every time instance on its forward counterpart. A robust treatment for
balancing the compute time versus trajectory storage has been
implemented in \PETSc{}, using 
the checkpointing algorithm {\em revolve}
\citep{griewank2000algorithm} as done by \citep{schanen_2016}. Previous work for ordinary differential equations (ODEs)
and
differential algebraic equations (DAEs) can be found in the work of \citep{zhang2014fatode}.

It is well known to the PDE-constrained optimization community that
continuous derivatives of the adjoint equation differ in nature from their discrete counterpart and may provide different gradients that
affect the convergence of the optimization algorithm \citep{amir2011reducing,limkilde2018reducing}. We do not
consider this issue in this paper; however, we note that the SUNDIALS package
\citep{sundials05} provides continuous adjoint capabilities
\citep{serban2005cvodes}. 

The paper is organized as follows. The time-dependent PDE-constrained
optimization problem is stated in an operator fashion that allows
for operators including diffusion, advection, and nonlinear ones such
as the $ u \cdot \nabla u $ operator in Burgers
equation. 
Then the main aspects of a spectral element method are presented, including a new approach for efficient matrix-free application of the transpose Jacobian matrix-vector product.
This is followed by a discussion of the discrete adjoint approach. We
then present an overview of the \PETSc{} ODE/DAE integrators and their adjoint
integrations, followed by a description of \TAO{} and its gradient-based
solvers. The results section focuses on the complete solution process for
performing PDE-constrained optimization using the spectral element
method and \PETSc{}/\TAO{}. We include a scalability study with over 2 billion unknowns and 130,000 cores. The final section briefly
summarizes our conclusions and potential future work.

\section{Problem formulation\label{sec:stateprob}}
The aim of this work is to illustrate a highly efficient spatial discretization for both the PDE and its adjoint
and, by connecting the various components available in the PETSc software stack, showcase the capability of performing large-scale 
PDE-constrained optimization. To this end, we present the mathematical framework at the same level of generality as handled by PETSc.

We denote a generic unsteady  PDE model by
\begin{eqnarray}
   \label{eq:generic}
 \mathbf{u}_t &=&  P[ \mathbf{u} ] , \  x \in \Omega \\ \nonumber
   B[\mathbf{u}|_{\partial \Omega}](t)&=& \mathbf{u}_b(t), \\ \nonumber
   \mathbf{u}( x,0)&=& \mathbf{u}_0( x),\nonumber
\end{eqnarray} 
where $ P[ \mathbf{u} ]$ is a stand-in operator for derivatives controlling the
spatial behavior of the solution $\mathbf{u}$ (this is not the most general notation since for some PDEs there may be also a dependence on
the spatial coordinates). The boundary condition $\mathbf{u}_b$ is
provided by an operator $B[\mathbf{u}]$, which is a scalar quantity for Dirichlet boundary conditions or a derivative for Neumann boundary
conditions. This notation does not exclude periodic boundary
conditions. For simplicity, in most of the presentation we assume
only homogeneous Dirichlet and/or periodic boundary conditions. However, this is not a limitation on the algorithms in \PETSc{}, since spatial discretizations are provided by the user.


Unless otherwise indicated, we  assume that all functions in this manuscript are in $\mathcal{K}$, a subspace of 
$L^2(\Omega)=\lbrace \mathbf{u} :\Omega\rightarrow \mathbb R |\quad (\int_{\Omega} |\mathbf{u}|^2 \d \Omega )^{1/2}<\infty\ \rbrace$. 
We assume that the subspace $\mathcal{K}$ is chosen such that the boundary conditions of the differential equations are satisfied.
One can express boundary conditions of the PDE as
additional constraints; however, in this work, we presume they are embedded in the discretization.

The optimization problem involves finding the control variable $\mathbf{p}$ that minimizes an objective functional ${\mathcal J}$ while satisfying the partial differential equation, Equation~\ref{eq:generic}. It can be stated as
\begin{eqnarray}
\label{eq:optimizationproblem}
\min_{p} &&{\mathcal J}(\mathbf{u},\mathbf{p}), \quad s.t. \\ \nonumber
\mathbf{u}_t&=&{P}[\mathbf{u}], \quad \text{with}\ \mathbf{u}|_{t=0}=\mathbf{u}_0 \in \mathcal{K}, \, \text{and} \  B[\mathbf{u}|_{\partial \Omega}](t)= \mathbf{u}_b(t),
\end{eqnarray}
where $\mathbf{u}$ is the state variable, which depends on $\mathbf{p}$.

In a general case an objective functional for time-dependent problems can be written as 
\begin{equation}
\mathcal J(\mathbf{u},\mathbf{p})=\int_0^T \int_{\Omega}g[\mathbf{u},\mathbf{p}] \delta(x-x_s, t-t_r)\d t \d \Omega\label{eq:costgeneral},
\end{equation}
where the Dirac delta restricts the objective to discrete points $x_s$ and $t_r$, at which an optimal solution is sought. 
We consider only PDE-constrained optimization problems where the control variables are the initial 
conditions $\mathbf{u}_0$, and we shall omit the $\mathbf{p}$ argument in the objective functional. However, general controls are fully supported by \PETSc{}/TAO{}.

A common occurrence
in the field of inverse problems is to have only a limited set of observations in space, at sensor locations 
$x_s$, in which case the objective functional defined above would take only a sparse set of values. One can seek solution trajectories that fit observations at discrete time instances, $t_r$, over a sparse set $R=\{t_r \in \mathbb R\}$. PETSc supports both scenarios; however, in this work, we consider a classical problem where $x_s$ is the entire set of discrete points describing the domain $\Omega$, and a single time observation.
If the objective functional is available only at time horizon $T$, then we can write
$$\mathcal J(\mathbf{u})=\int_0^T  \int_{\Omega}g[\mathbf{u}] \delta(\cdot, t-T)\d t \d \Omega =\int_{\Omega}g[\mathbf{u}(T)]\d \Omega\ .$$

An important special case is the data assimilation problem for
which we seek the initial condition $\mathbf{u}_0$ that leads at the {\em time horizon} $T$,  to a
solution $\mathbf{u}(T)$ that matches a reference solution $\mathbf{u}_d$. The standard objective 
functional that minimizes the difference between $ \mathbf{u} $ and the reference solution $ \mathbf{u}_d$  is
\begin{equation}
\label{eq:ourcost}  
  {\mathcal J}(\mathbf{u})=\int_{\Omega}(\mathbf{u}(T)- \mathbf{u}_d)^2 \ \d \Omega \,.  
\end{equation}

\section{Treatment of partial differential equations}
To outline the significant computational difficulties encountered in PDE-constrained optimization and also to showcase the full range of PETSc adjoint capabilities, we consider the viscous Burgers equation, a nonlinear problem in three dimensions.
For spatial coordinates $\mathbf{x}=[x_1,x_2,x_3]$ the periodic three-dimensional viscous Burgers equation can be stated as
\begin{align}
  \label{eq:burgers23d}
 \frac{\d \mathbf{u} }{\d t} &=\nu \Delta \mathbf{u} - \mathbf{u}\cdot \nabla \mathbf{u} \ ,
\end{align}
where $\nu$ represents the viscosity and $\mathbf{u}=[u_1, u_2, u_3]$. By casting the equation in the generic form of Equation~\ref{eq:generic} we have
\begin{align}
  \label{eq:burgers23dop}
 \frac{\d u_i }{\d t} &=\mu \Delta u_i-\sum_j( u_j\nabla_j u_i) \ , i,j=1,\ldots,3 \ ,
\end{align}
where $\nabla_j=\partial_{x_j}u_i$. The operator $
P[u_i]=\mu \Delta u_i-\sum_j( u_j\nabla_j u_i)$ incorporates all the problems we treat here and can be reduced either to a pure 
diffusion problem by removing the gradient term  
or to a linear advection problem by replacing the advection speed $u_i$ with a constant velocity $a_i$.

\subsection{Temporal discretization in PETSc}

The PETSc framework targets solutions adaptable to any time integration strategy.
It is, therefore, natural to consider a semi-discretization of the partial differential operators. 
Let us denote the spatial discretization of the operator $P[\mathbf{u}]$
as $\overline{P}[\overline{\mathbf{u}}]$, where $\overline{\mathbf{u}}$ is the
discretized $\mathbf{u}$ field.
Then the semi-discretization is 
\begin{equation}
\frac{\d \overline{\mathbf{u}}}{\d t} - \overline{P}[\overline{\mathbf{u}}] = 0 \ , \label{eq:semid}
\end{equation}
which is a form that can encapsulate both PDEs and ODEs. Although in the current work we treat only
 PDEs, we will at times refer to the semi-discrete form as an ODE.

Equation~\ref{eq:semid} is integrated in time either 
explicitly or implicitly by using a standard numerical integrator
\citep{tspaper}.
For implicit methods it is important to have access to the Jacobian (or its matrix-free application) $$J[\overline{\mathbf{u}}]= \frac{\d P[\overline{\mathbf{u}}]}{\d \overline{\mathbf{u}}}\ , $$ which can be derived from the semi-discrete form, as will be outlined in
Section~\ref{sec:sem}. The Jacobian is required only for implicit solvers or for the
discrete adjoint equation, as will be seen in Section~\ref{sec:opt}.

\subsection{Spatial discretization}
\label{sec:sem}
In this work, we focus on discretizations stemming from methods based on variational
(weak) formulations, such as finite elements or
spectral elements. This choice is justified by the flexibility and computational efficiency of such discretizations, considerations that are paramount to PDE-constrained optimization problems.

For the weak form of a partial differential equation we seek $ u $ in $\mathcal{K}$ with the property that 
the residual is orthogonal to the set of all test functions,
that is,
\[
\int_{\Omega} ( u_t- P[ u])  \ v \ \d \Omega  = 0
\]
for all $ v$ in $\mathcal{K}$. We illustrate here the prerequisites on the one dimensional operator $ P[u] = \nu \Delta u - u \cdot \nabla u$ and extend the discussion to higher dimensions in Section~\ref{sec:kron}. We apply integration by parts to obtain the continuous Galerkin formulation
\begin{equation}
\int_{\Omega} u_t v \ \d \Omega + \int_{\Omega}  v( u \cdot \nabla u)\ \d \Omega - \int_{\Omega}
\nu \nabla u \nabla v\ \d \Omega + \int_{\partial
  \Omega}  v \frac{\partial u}{\partial n}\  \d \partial \Omega =0,
\label{eq:weakcd}
\end{equation}
where $ n$ is the outward-facing normal.
The boundary term vanishes since $ u,  v \in \mathcal{K}$.

Based on the weak form, several
discretizations are suitable, including the finite element method, spectral
element method,  or discontinuous Galerkin method.  The spectral
element method is a subclass of Galerkin methods, or weighted residual
methods, that minimize the error of
the numerical solution in the energy norm over a chosen space of
polynomials or, equivalently, require the error to be orthogonal to
the subspace defined by the spectral elements. These concepts are discussed in detail in \citep{fischer:hom}.

In the following we illustrate the spectral element discretization for  Equation~\ref{eq:weakcd}.
The domain $\Omega=\cup_{e=1,M} \Omega_e $ is decomposed into $M$ non-overlapping 
subdomains $\Omega_e$, termed elements, over which the data will be represented by orthogonal polynomials.

The space of polynomials of order $N$ defined over an element $\Omega_e, \ e=1,\ldots, E$, is
$$
\mathbb P_{N,E}=\lbrace  \phi| \phi \in L^2(\Omega); \quad \phi|_{\Omega_e} \text{polynomial of degree} \leq N\rbrace.
$$
Subsequently $\mathcal{P}_N=\mathcal{K}\cap  P_{N,E}$.
In the polynomial space $\mathcal{P}_N$ we represent the numerical solution as $u( x)=\sum_{i=0}^N u_{i}\phi_i( x)$.

The discrete point distribution $ x$ can be a Chebyshev grid
or a Legendre grid, consistent with the polynomial discretization. The current work is performed using Legendre polynomials and Gauss-Legendre-Lobatto (GLL) grids, since this representation imposes less rigid stability restrictions
on the time-steppers, having a smaller clustering of grid points at the boundary ends.

To proceed with the numerical discretization, we plug the ansatz into
Equation~\ref{eq:weakcd} to obtain
\begin{multline}
  \frac{\partial}{\partial t}\sum_{i,j=0}^N u_{i}v_{j}\underbrace{\int_{\Omega} \phi_i(x)\phi_j(x) \d \Omega}_{M_{ij}} - 
  \nu\sum_{i,j=1}^N u_{i} v_{j}\underbrace{\int_{\Omega} \phi'_i(x) \phi'_j(x) \d \Omega}_{K_{ij}} +\\
  \sum_{i,j,k=1}^N u_{k}u_{i} v_{j}\underbrace{\int_{\Omega}\phi_k(x)\phi'_i(x) \phi_j(x) \d \Omega}_{D_{ij}}  =0,
\end{multline}
where $M_{ij}$ is the mass matrix, which is diagonal in this case because of the orthogonality property of the polynomials,
$K_{ij}$ is the stiffness
matrix, and $D_{ij}$ is the differential matrix.
Note that in this case, the contribution of the basis $\phi_k$ is
present only as integration weights since the polynomials $\phi_k$ evaluate to unity on the $\Omega$ grid. Given that the stiffness matrix can be obtained from the differential as $K=D^\top M D,$ various implementations are possible, according to the complexity of the physical problem or geometry. For example, the CEED library \citep{ceed_report} decomposes $K$ in all situations and has a highly performant implementation on heterogeneous architectures while, for speed on rectangular geometries, Nek5000 \citep{nekdoc} has options which implement the forward operator $K$ as discussed in this work.

We can write the equations in algebraic form and scale out the test function $v$ to obtain 
\begin{equation}
  \label{eq:discrete1d}
  M\frac{d \overline{u}}{d t} =-\nu\ K \overline{u}+\overline{u} \circ D \overline{u},
\end{equation}
where $\overline{u}=(u_0,\ u_1, \ldots, u_N)$.
Here we introduce the notation $\circ$ to designate the Hadamard product, also known as pointwise multiplication.

If we assume  $\Omega_e = [a, \ b]$ and define the reference element $\hat{\Omega} = [-1, \ 1]$,
then for $x\in \Omega_e$ and $r \in \hat{\Omega}$ the mapping from each element to the reference element is
$$x= a+ \frac{(b-a)(r+1)}{2}\ .$$
This mapping introduces a scaling factor for each element, given by the Jacobian of the coordinate transformation. For curvilinear geometries, as well as for variable coefficients $\nu$, the operators $K$ and $D$  have more complex expressions that take into account the Jacobian and multiply the variable coefficients. Since our goal is to outline a vectorized implementation of the adjoint equations, we focus on simple affine transformations to facilitate the presentation. As a result, for a domain $\Omega=[0, \ L_1]$, the operators $M$ and $K$ gain a simple scaling, while the operator $D$ is unchanged by affine transformations, and we  introduce the following notation
\begin{eqnarray*}
M_1&=&\frac{L_1}{2}M\\
K_1&=&\frac{2}{L_1}K.
\end{eqnarray*}

For an implicit integrator, Equation~\ref{eq:discrete1d} requires the Jacobian of the right-hand side. The Jacobian is required in its transposed form by the discrete adjoint 
integrator.
For the right-hand side $\overline{P}[\overline u]= \nu\ K \overline{u}- \overline{u} \circ D \overline{u}$, the Jacobian becomes
\begin{equation}
\label{eq:jac1d}
J=\frac{d \overline P}{d \overline u}=\nu \ K - (\diag(D \overline u) + \diag(\overline u) D),
\end{equation}
where we use  the derivative of a Hadamard product to compute the derivative of the nonlinear term.
Consider the Jacobian applied to a field $w$,
\begin{eqnarray}
\label{eq:jac1dmatfree}
Jw&=&(\nu \ K - (\diag(D [\overline u]) + \diag(\overline u) D))[w]\\ \nonumber
&=&\nu \ K[w] - w\circ D [\overline u] + \overline u \circ D[w]\,. \nonumber
\end{eqnarray}

The Jacobian transpose that will generate the adjoint in the optimization problem is given in matrix-free form as
\begin{eqnarray}
\label{eq:jactrans1dmatfree}
J^\top w&=&\nu \ K^\top[w] - w\circ D [\overline u] + \overline u \circ D^\top[w]\,. \nonumber
\end{eqnarray}
This particular matrix-free approach brings impressive benefits in three dimensions.

\subsection{Vectorized matrix-free discretization}\label{sec:kron}
Consider the domain $\Omega=[-1,\ 1]^d$ discretized in $N$ GLL points, $d=3$. The  basis function is
 separable $\phi_k(\mathbf{x})=\phi_i(x_1)\phi_j(x_2)\phi_l(x_3),$ where $i,j,l=1,\ldots,N$.
The ansatz on the solution is like the one-dimensional case:
$$\overline{\mathbf{u}}(\mathbf{x})=\sum_{i,j,l=0}^N u_{ijl}\phi_i(x_1)\phi_j(x_2)\phi_l(x_3).$$

The separability of the basis functions allows us to represent the operators $\mathbf{K}, \ \mathbf{D},$ and $ \mathbf{M}$ as tensor products in each dimension, as will be shown, this brings high computational advantages. Consider a domain $\Omega=[0, L_1] \times [0, L_2]\times [0, L_3],$ and consider a second derivative operator applied to a component $u_i$, for example $\partial_{yy}\overline u_i\approx \mathbf{K}_2[ \overline u_i]=(M_1\otimes K_2\otimes M_3)\overline u_i$. Note that we use boldface operators for vector fields, i.e. $\mathbf {K}_i$ applies the second order derivative to a single component $\overline u_i(\mathbf{x})=u_i(x_1,x_2,x_3)$ of the three-dimensional field $\overline{\mathbf{u}}=(u_1,u_2,u_3)$, that is, $K_j$ applies only to direction $x_j$. Summing all components, we obtain the full Laplacian in Equation~\ref{eq:burgers23d} as applied to component $u_i$, 

$$\Delta u_i=(\partial_{xx}+\partial_{yy}+\partial_{zz})[u_i]\approx (\mathbf{K}_1 + \mathbf{K}_2+\mathbf{K}_3)][\overline u_i], $$ which can be implemented by using
\begin{align*}
  \mathbf K[\overline u_i] =[(K_1\otimes M_2\otimes M_3) + (M_1\otimes K_2\otimes M_3)+(M_1\otimes M_2\otimes K_3)][\overline u_i].
\end{align*}
The advantage of this representation is that the operators can be evaluated efficiently. Each matrix-vector multiplication, without unrolling the unknown 
variable $u_i$, is  replaced  with matrix-matrix products, for example, $\mathbf K_2[\overline u_i]=(M_1 \otimes K_2\otimes M_3) \overline{u_i} = M_3^\top   \overline{u_i} K_2M_1$.
We can now state the discrete form of Equation~\ref{eq:burgers23d} with $P[\overline u_i]=\mu \Delta\overline u_i-\sum_j(\overline u_j\nabla_j \overline u_i)$ as
\begin{align}
  \label{eq:discrete3d}
 \mathbf M \frac{\d \overline{ u}_i}{\d t} &=\underbrace{\nu\ \mathbf K[\overline{ u}_i]}_{\mathcal L (\overline{u}_i)}- 
 \underbrace{\overline{\mathbf u}\cdot\mathbf G[\overline{u}_i]}_{\mathcal N (\overline{u}_i)},
\end{align}
where $\mathbf G=[\mathbf D_1,\mathbf D_2, \mathbf D_3]$, and the operators are given by
\begin{align}
\mathbf M &=M_1\otimes M_2\otimes M_3 \\ \nonumber
\mathbf K &=(K_1\otimes M_2\otimes M_3) + (M_1\otimes K_2\otimes M_3)+(M_1\otimes M_2\otimes K_3) \\ \nonumber
\mathbf G &=[(D_1\otimes M_2\otimes M_3), (M_1\otimes D_2\otimes M_3),(M_1\otimes M_2\otimes D_3)]. \nonumber
\end{align}

The discrete right-hand side $\overline{P}[\overline u_i]$ can be implemented efficiently by using the Kronecker product evaluations as
\begin{equation}\label{eq:rhsvec}
\overline P[\overline u_i]=\sum_j \mathbf K_j[\overline u_i]+\sum_j\overline u_j\circ \mathbf D_j [\overline u_i],
\end{equation}
where each operator $\mathbf K_j[\overline u_i]$ and $\mathbf D_j [\overline u_i]$ is a tensor contraction implemented in a vectorized fashion by using BLAS calls.

To take advantage of these efficient evaluations, we need to implement the Jacobian and Jacobian transpose also in matrix-free form. Let us first analyze the Jacobian matrix, which for the system in Equation~\ref{eq:discrete3d} is given as 
$$
J_{il}=\frac{\d\overline P[\overline u_i]}{\d\overline u_l},\ i,l=1,\ldots, 3\,,
$$
where each component $i,l$ is a block matrix, symmetric for the first component, $\mathcal L (\overline{u}_i)$, corresponding to the Laplacian operator $\mathbf{K}$, but asymmetric for the components corresponding to the nonlinear operator $\mathcal N (\overline{u}_i)$. The Jacobian applied to a vector $\overline{\mathbf w}=(\overline w_1, \ \overline w_2, \ \overline w_3)$ has the representation
\begin{equation*}
J \overline{\mathbf w}=   \begin{bmatrix}\mathbf K +\frac{d \mathcal N(\overline u_1)}{d \overline u_1}& \frac{d \mathcal N(\overline u_1)}{d\overline u_2} & \frac{d \mathcal N(\overline u_1)}{d\overline u_3}\\
\frac{d \mathcal N(\overline u_2)}{d\overline u_1}&  \mathbf K +\frac{d \mathcal N(\overline u_2)}{d\overline u_2}&\frac{d \mathcal N(\overline u_2)}{d\overline u_3}\\
\frac{d \mathcal N(\overline u_3)}{d\overline u_1}&  \frac{d \mathcal N(\overline u_3)}{d\overline u_2}& \mathbf K + \frac{d \mathcal N(\overline u_3)}{d\overline u_3}
\end{bmatrix}
\begin{bmatrix} 
w_1 \\
w_2 \\
w_3
\end{bmatrix}.
\end{equation*}

For each derivative term in the Jacobian we apply the chain rule of matrix calculus to obtain
\begin{eqnarray*}
J_{il}=\frac{\d\overline P[\overline u_i]}{\d\overline u_l}&=&\frac{\d}{\d\overline u_l}\sum_j \mathbf  K_j[\overline u_i]+\frac{\d}{\d\overline u_l}\sum_j\overline u_j\mathbf  D_j[\overline u_i]\\
&=&\mathbf  K_j[\delta_{il}]+\sum_j\left(\frac{\d}{\d\overline u_l}\overline u_j\mathbf  D_j[\overline u_i] + \overline u_j\frac{\d}{\d\overline u_l}\mathbf  D_j[\overline u_i] \right)\\
&=&\mathbf  K_j[\delta_{il}]+\sum_j\delta_{jl}\mathbf  D_j[\overline u_i] + \sum_j\overline u_j\mathbf D_j[\delta_{il}]\ , i,l=1,\ldots,3, 
\end{eqnarray*}
where we use the fact that $\frac{\d\overline u_i}{\d\overline u_l}=\delta_{jl}$, with $\delta_{jl}$ the Kronecker delta. By applying $J$ to a vector field $\overline w_l$, with $l=1,\ldots,3$, we obtain
\begin{eqnarray*}
J\overline{\mathbf w}&=&\sum_l\sum_j \mathbf K_j[\delta_{il}\overline w_l]+\sum_l\sum_j\mathbf  D_j[\overline u_i]\delta_{jl}\overline w_l + \sum_l\sum_j(\overline u_j\mathbf D_j[\delta_{il}])\overline w_l\\
&=&\sum_j\mathbf K_j[\overline w_i]+\sum_j(\mathbf  D_j[\overline u_i]\overline w_j+  \overline u_j\mathbf D_j[\overline w_i]).
\end{eqnarray*}

For consistency with Equation~\ref{eq:rhsvec} we state the vectorized matrix-free form as
\begin{equation}\label{eq:jacvec}
J\overline{\mathbf w}=\sum_j \mathbf K_j[\overline w_i]+\sum_j(\overline w_j\circ \mathbf D_j [\overline u_i]+\overline u_j\circ \mathbf D_j [\overline w_i]).
\end{equation}

To compute the transpose, we need not only consider $J_{il}\rightarrow J_{li}$ but also transpose blocks $\mathbf D_i$, which are asymmetric in nature, yielding
\begin{eqnarray*}
J^\top =J_{li}&=&\sum_j \mathbf K_j^\top[\delta_{li}]+\mathbf  D_l[\overline u_i] + \sum_j \mathbf D_j^\top \delta_{li}\overline u_j^\top,
\end{eqnarray*}
where we used $\sum_j\mathbf  D_j[\overline u_i]\delta_{jl}=\mathbf  D_l[\overline u_i]$. Applying this result to $\mathbf w$ and using Einstein summation, as well as consolidation of indices, with $i$ transformed into $j,$ we have
\begin{eqnarray*}
J^\top \overline{\mathbf w}&=&\sum_i \sum_j \mathbf K_j^\top[\delta_{li}\overline w_i]+\sum_i \mathbf  D_l[\overline u_i] \overline w_i + \sum_i \sum_j( \mathbf D_j^\top \delta_{li}\overline u_j^\top)\overline w_i\\
&=&\sum_j \mathbf K_j^\top[\overline w_l]+\sum_j\mathbf  D_l[\overline u_j] \overline w_j+ \sum_j \mathbf D_j^\top \overline w_l^\top \overline u_j^\top.
\end{eqnarray*}
We conclude with a compact form for the matrix-free representation of the Jacobian transpose:
\begin{equation}\label{eq:jactransvec}
J^\top\overline{\mathbf w}=\sum_j \mathbf K_j^\top[\overline w_i]+\sum_j(\overline w_j\circ \mathbf D_i [\overline u_j]+\mathbf D_j^\top [\overline u_j\circ \overline w_i]).
\end{equation}
The computational complexity of Kronecker product evaluations lowers the cost of simple block evaluations from $\mathcal O (N^{2d})$ to  
$\mathcal O (N^{d+1})$, where $d$ is the spatial dimension, thus having an increasing impact with increasing dimensionality. This can be exploited further by considering tensorizations over elements of a mesh, as done by \citep{ceed_report}, or additional fields, such as pressure and temperature in compressible Navier-Stokes.

In the discrete optimization community, the Jacobian, Equation~\ref{eq:jacvec}, is usually referred to as the forward mode, while the transpose of the Jacobian, Equation~\ref{eq:jactransvec}, is called the reverse mode. While the computation of the forward mode allows for several methods, the reverse mode is often computed by using an automatic differentiation (AD) strategy, as provided by libraries such as ADOL-C \citep{adolc}, or CoDiPack \citep{codipack}. AD computes the derivatives by applying the chain rule at the instruction level and requires time for the reverse mode of at least twice the time for the forward mode. The derivation of the Jacobian transpose we provide here is straight forward to implement and preserves the vectorization of the forward integrator. In Section 6 we that demonstrate the time can be less than twice the time of the forward integration. This suggests that AD libraries applied to PDEs may benefit from applying the chain rule, as described in the current work, at the matrix block level, that is, the chain rule of matrix calculus, rather than at the instruction level.

\subsection{PETSc implementation}

The \PETSc{} time-stepping component \TS{} provides access to numerous ODE
integrators including explicit, implicit, and implicit-explicit
methods. 
%
The linear systems that arise in the
implicit solvers may be solved by using any of \PETSc{}'s algebraic solvers,
as well as solvers from other packages such as hypre \citep{hypre-users-manual} or SuperLU\_Dist \citep{lidemmel03}. The \TS{}
integrators in PETSc are generally multistage integrators, with local and
global error estimate adaptive controllers, although \PETSc{} does also
provide some multistep integrators. Details on solvers and integrators in \PETSc{} may be found in the user's manual \citep{petsc-user-ref}.

The \PETSc{} interface for solving time-dependent problems is organized
 around the following form of an implicit differential equation:
\[
        F(t,u,\dot{u}) = G(t,u), \quad u(t_0) = u_0.
\]
If the matrix $F_{\dot{u}}(t) = \partial F / \partial \dot{u}$ is
nonsingular, then this is an ordinary differential equation and can be transformed to the standard
explicit form. For PDEs we can often write $F(t,u,\dot{u})=M \dot{\overline
  u}$ and $G(t,u)=\overline P[\overline u]$, or  $F(t,u,\dot{u})=\dot{\overline
  u}$ and $G(t,u)=M^{-1}\overline P[\overline u]$, since the representation 
of any semi-discrete partial differential equation is 
an ODE. 
 For ODEs, with nontrivial mass
matrices, the
above implicit/DAE interface can significantly reduce the overhead in preparing
the system for algebraic solvers by having the user code
assemble the correctly shifted matrix instead of having the library do it with less efficient general code.  The user provides function pointers
and pointers to user-defined data for each operation, such as
the function $ F()$ needed by the library. This approach allows full utilization from C,
Fortran, Python, and C++ since all these languages support these
constructs, whereas, for example, the use of C++ classes would limit the
language portability.

The PETSc TS object manages all the time-steppers in PETSc. By selecting various options, it supplies all the combinations of solvers. The PETSc nonlinear solver object supplies the user interface for nonlinear systems that, using Newton-type methods, rely on linear system solves handled by the PETSc Krylov solver object. The Jacobian and Jacobian transpose operators can be provided by utilizing sparse matrix storage, for example, compressed sparse row storage, or matrix-free implementations. Although traditional uses of PETSc target linear solvers based on sparsely stored matrices, the matrix-free approach can be used for both Jacobian and its transpose and is chosen in the current context to attain higher efficiency of the adjoint solution computation. In PETSc, a matrix provided as the action of an operator is represented as a  \texttt{Shell} matrix, on which an operation as defined as below
\begin{lstlisting}
MatCreateShell(MPI_Comm comm,PetscInt m,PetscInt n,PetscInt M,PetscInt N,void *ctx,Mat *A);
MatShellSetOperation(Mat A,MATOP_MUL, (*g)(Mat,Vec,Vec))
\end{lstlisting}

Time integrators in PETSc have a common interface for both implicit and explicit methods and combinations thereof, as in the following example listing.
\begin{lstlisting}
TSCreate(MPI_Comm comm, TS *ts); 
TSSetRHSFunction(TS ts,Vec r,(*f)(TS,PetscReal,Vec,Vec,void*),void *ctx);
TSSetRHSJacobian(TS ts,Mat mat,Mat pmat, (*j)(TS,PetscReal,Vec,Mat,Mat,void*),void *ctx);
TSSolve(TS ts,Vec u);
\end{lstlisting}
The matrix object is passed to the time integrator with \texttt{TSSetRHSJacobian}.

\section{PDE-constrained optimization}

Several mathematical approaches are available in PDE-constrained
optimization, of which the Lagrange-multiplier
framework can be used in both continuous and discrete adjoint-based optimization.

\subsection{Using adjoints for optimization}
Equation~\ref{eq:optimizationproblem} is posed as a minimization of a functional subject to a time-dependent PDE-constraint.
For continuous adjoints, the appropriate inner product
space is defined by
\begin{equation}
(\mathbf{u}\ ,\mathbf{v})=\int_0^T  \int_{\Omega}\mathbf{u} \cdot \mathbf{v} \ \d \Omega \d t .
\label{eq:inner}
\end{equation}

 Equation~\ref{eq:optimizationproblem} can now be written in the Lagrangian framework, with $\mathbf{v}$
playing the role of Lagrange multiplier or adjoint variable:
\begin{subequations}
\begin{align}
\label{unconstrainedoptimizationproblem}
  \Lagrangian(\mathbf{u}_0,\mathbf{u},\mathbf{v}) = \mathcal J(\mathbf{u}) + 
  \int_0^T\int_{\Omega}(\mathbf{v} (\mathbf{u}_t- P[\mathbf{u}])) \d \Omega \d t\,.
\end{align}
\end{subequations}

The Lagrangian has a stationary point, $\nabla \mathcal L=0$, for each extremum of the original problem and provides a necessary condition for optimality of solution (additional stationary points may exist that are not solutions to the optimization problem).

For a semi-discrete partial differential equation $\overline u_t = \overline{P}[\overline u]$
let us consider a general model of an explicit time-stepper:
\begin{equation}
\overline{\mathbf{u}}_{n+1}=A(\overline{\mathbf{u}}_n,\overline P[\overline{\mathbf{u}}_n]).
\end{equation}
For implicit methods we need
to consider the right-hand side as an operator applied to $\overline{\mathbf{u}}_{n+1}$. 
In the discrete framework, the inner product defined in Equation~\ref{eq:inner} has to be
discretized consistently compatible with the discretization of the PDE. For a Galerkin method, the discrete inner product in space is $(u,v)=u^\top  M v$, where $M$ is the mass matrix, and for
the time integration we utilize the Riemann sum. The total inner product in 
Problem~\ref{unconstrainedoptimizationproblem} is
\begin{equation*}
\int_0^T\int_{\Omega}(\mathbf{v}(\mathbf{u}_t- P[\mathbf{u}])) \d \Omega \d t \approx \sum_{n=1}^N\mathbf{v}_n^\top  M (\overline{\mathbf{u}}_n-
 A(\overline{\mathbf{u}}_{n-1},\overline P[\overline{\mathbf{u}}_{n-1}])),
\end{equation*}
where we use the convention that $\mathbf{u}_0=\mathbf{u}(0)$ and $\mathbf{u}_N =\mathbf{u}(T)$.

After a shift and rewrite of the summation bounds, the total Lagrangian is
\begin{eqnarray}
  \Lagrangian(\mathbf{u}_0,\overline{\mathbf{u}},\overline{\mathbf{v}}) = \overline{\mathcal J}(\overline{\mathbf{u}}) + \sum_{n=1}^{N}\overline{\mathbf{v}}_n^\top  M \overline{\mathbf{u}}_n-\sum_{n=0}^{N-1}
 \overline{\mathbf{v}}_{n+1}^\top  M A(\overline{\mathbf{u}}_n,\overline P[\overline{\mathbf{u}}_n]).
\end{eqnarray}

To identify the adjoint, we now require that all derivatives of $\mathcal L$ cancel,
\begin{align}
 \label{eq:adjointloop}
  \frac{\partial \Lagrangian}{\partial \overline{\mathbf{u}}_n} &= \frac{\partial \mathcal J}{\partial \overline{\mathbf{u}}_n} -\overline{\mathbf{v}}_n^\top  M
  +\bigg(\frac{\partial A(\overline{\mathbf{u}}_n,\overline P[\overline{\mathbf{u}}_n])}{\partial \overline{\mathbf{u}}_n}\bigg)^\top  \overline{\mathbf{v}}_{n+1}^\top  M=0\\
  \label{eq:initadj}
   \frac{\partial \Lagrangian}{\partial \overline{\mathbf{u}}_N} &= \frac{\partial \mathcal J}{\partial \overline{\mathbf{u}}_N} -\overline{\mathbf{v}}_N^\top  M =0\\
    \frac{\partial \Lagrangian}{\partial \overline{\mathbf{u}}_0} &= \frac{\partial \mathcal J}{\partial \overline{\mathbf{u}}_0} 
  +\bigg(\frac{\partial A(\overline{\mathbf{u}}_n,\overline P[\overline{\mathbf{u}}_0])}{\partial \overline{\mathbf{u}}_0}\bigg)^\top  \overline{\mathbf{v}}_{1}^\top  M=0.
  \label{eq:gradient}
\end{align}

The adjoint equation to be solved is provided by  Equation~\ref{eq:adjointloop}, while Equation~\ref{eq:initadj} gives the 
initial condition for the adjoint equation and Equation~\ref{eq:gradient} is the gradient to be used in the optimization step.

We now can analyze  the operator $\frac{\partial A(\overline{\mathbf{u}},\overline P[\overline{\mathbf{u}}])}
{\partial \overline{\mathbf{u}}}$. To take the derivative of the operator $A$, we treat it via the
implicit function theorem where $\overline{\mathbf{u}}$ and $\overline P[\overline{\mathbf{u}}]$ are regarded as independent variables yielding
\begin{equation}
\label{eq:adjtime}
\bigg(\frac{{\rm d}A(\overline{\mathbf{u}},\overline P[\overline{\mathbf{u}}])}{{\rm d} \overline{\mathbf{u}}}\bigg)^\top =
\bigg(\frac{\partial A(\overline{\mathbf{u}},\overline P[\overline{\mathbf{u}}])}{\partial \overline{\mathbf{u}}}+ \frac{\partial A(\overline{\mathbf{u}},\overline P[\overline{\mathbf{u}}])}{\partial \overline P[\overline{\mathbf{u}}]}\frac{\partial \overline P[\overline{\mathbf{u}}])}{\partial \overline{\mathbf{u}}} \bigg)^\top ,
\end{equation}
where, for example, for one-dimensional spatial operators, such as $P[\overline u]= K [\overline u] + u \circ D [\overline u ]$, the Jacobian is given by 
Equation~\ref{eq:jac1d}, while the adjoint formulation relies on the Jacobian transpose of Equation~\ref{eq:jactrans1dmatfree}. For three-dimensions the matrix-free Jacobian and Jacobian transpose are computed via Equation~\ref{eq:jacvec}, and Equation~\ref{eq:jactransvec}, respectively.

The Jacobian of the operator $A$ depends on the time integration scheme and leads to the tableau of the time integrator. We illustrate this on the 
forward Euler time integration scheme. The model $\overline{\mathbf u}_{n+1}=A(\overline{\mathbf u}_n,\overline P[\overline{\mathbf u}_n])$ becomes
\begin{equation*}
\overline{\mathbf u}_{n+1}= \overline{\mathbf u}_n+\Delta t \overline P[\overline{\mathbf u}_n].
\end{equation*}
The derivative of $A$ is calculated by using Equation~\ref{eq:adjtime} and becomes
\begin{equation}
\frac{\partial A(\overline{\mathbf u},\overline P[\overline{\mathbf u}])}{\partial \overline{\mathbf u}}=
I+ \Delta t \frac{\partial \overline P[\overline{\mathbf u}]}{\partial \overline{\mathbf u}}.
\end{equation}
Setting $\nabla_{u_n} \mathcal L=0$ in Equation~\ref{eq:adjointloop}, we 
obtain the adjoint equation based on the forward Euler with the time reversal $n=N,\ldots,1$:
\begin{equation}
\overline{\mathbf v}_n= \overline{\mathbf v}_{n+1}+\Delta t \frac{\partial \overline P[\overline{\mathbf u}]}{\partial \overline{\mathbf u}_n}^\top  \overline{\mathbf v}_{n+1}.
\end{equation}
For more details about the discrete adjoint implementation in PETSc, see the
derivations in \citep{Zhang_2019}.

\begin{figure}[!htb]
  \begin{center}
    \begin{tikzpicture}[scale=1.8]
      \node (A) at (0,1) {$\overline{\mathbf u}_0$};
      \node (a) at (2,1) {$\overline{\mathbf u}_k$};
      \node (b) at (4,1) {$\overline{\mathbf u}_n$};
      \node (B) at (6,1) {$\overline{\mathbf u}_N$};
      \node (C) at (6,0) {$\overline{\mathbf v}_N (\nabla_{\overline{\mathbf u}_N}  \mathcal L=0)$};
      \node (c) at (4,0) {$\overline{\mathbf v}_{n}$};
      \node (d) at (2,0) {$\overline{\mathbf v}_{k}$};
      \node (D) at (0,0) {$\overline{\mathbf v}_0 (\nabla_{\overline{\mathbf u}_0}  \mathcal L=0)$};
      \draw [->,dashed] (A) edge (a) (a) edge (b) (b) edge (B);
      \draw [->] (B) edge (C);
      \draw [->,dashed] (C) edge (c) (c) edge (d) (d) edge (D);
      \draw [->] (b) edge (c) (a) edge (d);
      \draw [->] (D) edge (A);
    \end{tikzpicture}
  \end{center}
  \caption{Workflow for optimization problems using adjoints. 
    \label{fig:algorithm}}
\end{figure}
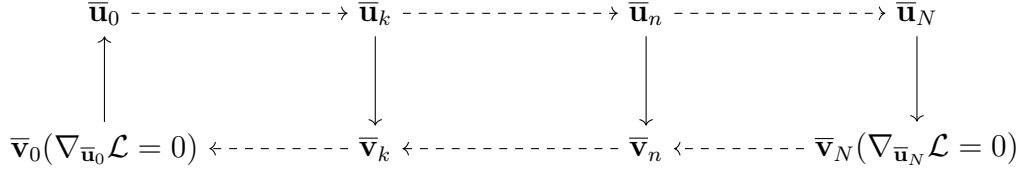

This procedure requires that the
forward solution be available during the backward solve. This data management procedure is referred to as checkpointing and 
is available directly via the PETSc adjoint implementation, as will be described in Section~\ref{sec:checkpointing}.
We summarize the adjoint-based approach in Figure~\ref{fig:algorithm}, where
initial conditions for the adjoint integration are provided by
the solution of Equation~\ref{eq:initadj}.

\section{Computational considerations}

\subsection{Accuracy}

The time integration error in \PETSc{} can be tracked by using the built-in error
estimation and error control mechanism.  
These work by changing the step size to maintain user-specified absolute, $\TolA $, and relative, $\TolR$,
tolerances. The error estimate is based on the local truncation error, so for
every step the algorithm verifies that the estimated local truncation
error satisfies the tolerances provided by the user and computes a new
step size. For multistage methods, the local truncation is
obtained by comparing the solution $\overline{\mathbf u}$ to a lower-order
$\widehat{p}=p-1$ approximation, $\widehat{\overline{\mathbf u}}$, where $p$ 
is the order of the method and $\widehat{p}$ the order of $\widehat{\overline{\mathbf u}}$.

The adaptive controller at step $n$ computes a tolerance level,
\begin{eqnarray*}
\Tol_n(i)&=&\TolA(i) +  \max(|\overline{\mathbf{u}}_n(i)|,|\widehat{\overline{\mathbf u}}_n(i)|) \TolR(i)\,,
\end{eqnarray*}
and forms the acceptable error level, referred to as the weighted local truncation error
\begin{align*}
\SErr_n= \frac{1}{m} \sum_{i=1}^{m}\sqrt{\frac{\left\|\overline{\mathbf{u}}_n(i)
  -\widehat{\overline{\mathbf u}}_n(i)\right\|}{\Tol_n(i)}} ~ \textnormal{or} ~ \SErr_n= \max_{1\dots m}\frac{\left\|\overline{\mathbf u}_n(i)
  -\widehat{\overline{\mathbf u}}_n(i)\right\|}{\Tol(i)}\,,
\end{align*}
where the errors are computed component-wise and $m$ is the dimension of
$u$. 
The error tolerances are satisfied when $\SErr \le 1.0$. 
The next step size is based on this error estimate and is determined by
\begin{eqnarray}
  \label{eq:hnew}
 \Delta t_{\rm new}(t)&=&\Delta t_{\rm{old}} \min(\alpha_{\max},
 \max(\alpha_{\min}, \beta (1/\SErr)^\frac{1}{\widehat{p}+1}))\,,
\end{eqnarray}
where $\alpha_{\min} $ and 
$\alpha_{\max}$ keep the change in
$\Delta t$
to within a certain factor and $\beta<1$ is chosen
so that there is
some margin which the tolerances are satisfied and so that the
probability of rejection is decreased. \PETSc{} also provides a global
error or a posteriori error estimation and control mechanism similar
to the local error control implementation \citep{constantinescu2015estimating}.

To assess the error incurred in the spatial spectral element discretization, we use the
a posteriori spectral error analysis suggested by
\citep{Mavriplis1990}. An exhaustive review is available
in \citep{kopriva2009}.
Considering the discrete solution $\overline{\mathbf{u}}\in \mathcal{P}_N$ (the space of
polynomials of order $N$), an
extrapolated solution $\tilde{\mathbf{u}}\in \mathcal{P}_M$, where $ {\mathcal{P}_m} $ is a  space of polynomials of
order $M\gg N$, and the exact solution $\mathbf{u}$, the error can be bounded as 
$$\Vert \mathbf{u} - \mathbf{u}_N\Vert \leq \underbrace{\Vert \mathbf{u}-\tilde{\mathbf{u}}\Vert}_{\text{extrapolation error}} + \underbrace{\Vert \tilde{\mathbf{u}}-
  \Pi_N \tilde{\mathbf{u}}\Vert}_{\text{truncation error}} + \underbrace{\Vert
  \mathbf{u}_N-\Pi_N \tilde{\mathbf{u}}\Vert}_{\text{quadrature error}} \ ,$$ where $\Pi_N
\mathbf{u}$ is the projection of $\mathbf{u}$ on a polynomial space of order $N$.

The \textit{extrapolation error} is the error incurred by
  extrapolating the solution to a higher-order space and determines the leading-order error. The \textit{truncation error}, due to the truncation to a lower-order space, and the \textit{quadrature error}, introduced by discretizing the weak formulation, are of a similar higher order. The leading-order term of the truncation in the error bound to $\epsilon_e$ is 
$$\epsilon_e=\frac{a_N^2}{(2N+2)/2}\ , $$ where $a_N$ can be
approximated from a least-squares expansion as $a_N=c\mathrm e^{-\sigma
  N}$, $N$ is the number of degrees of freedom, $\sigma$ is the slope of the decay, and $c$ is a scaling factor. The exponential decay of the error indicates that the same number of discretization points as for a low-order method yields far higher accuracy. Note that spectral accuracy is exhibited only within an element of a mesh, and across elements, continuity can only be assured in a continuous Galerkin framework. The exponential decay, combined with the fact that spectral element methods are amenable to highly vectorized implementations, renders them to be considered state of the art in terms of efficiency. Nevertheless, the spectral element method is not free of other drawbacks, such as poor representations of shocks and sharp interfaces and prohibitive stability restrictions.

\subsection{Computational efficiency}
Both the implementation and algorithm can have a high impact on the total time of the simulation.
The best approach is to choose a scheme that has the lowest number of 
flops per accuracy threshold, and these methods are the ones with spectral accuracy.
Of the spectral methods, the spectral element method has the significant advantage is that it is 
highly parallelizable and can easily perform a matrix-free
tensor product per spectral element.
The computational complexity of the spectral element method can be reduced from $\mathcal O (N^{2d})$ to  
$\mathcal O (N^{d+1})$, where $d$ is the spatial dimension. The highly vectorized structure of the tensor product multiplications
allows efficient implementations, see \citep{fischer:hom} p 168. Here we demonstrate a similar efficient approach for the
application of the transpose Jacobian with the same work estimates; see Section \ref{sec:kron}.

\subsection{Data management}
\label{sec:checkpointing}
To integrate the adjoint equation, one needs to solve Equation~\ref{eq:adjointloop} at every time step.   As
illustrated in Figure~\ref{fig:algorithm}, at every
timestep of the adjoint integration, the corresponding forward solution is
needed. Writing the solution to disk may be expensive, depending on the I/O implementation of the code and the machine specifications. Storing the solution in memory may also
be prohibitive depending on the problem size or the available memory.

These issues raise the question of whether it is more efficient to store
the solution of the forward problem to file or in memory or to recompute
the solution. Ideally, one would seek to find the break-even point between retrieving and recomputing the solution.
Sometimes, it may be advantageous to recompute the solution from a stored solution partially. The optimal distribution of recomputations per solution store, when the total number of time steps and storage space is known
beforehand, has been established to be given by the revolve algorithm \citep{griewank2000algorithm}.

\subsection{PETSc/TAO implementation}

The \PETSc{} interfaces for computing discrete adjoints are built on top of the ODE/DAE interface. We illustrate the
adjoint derivation only for initial conditions; however, the framework allows also for optimization with respect to other
parameters.
The adjoint integrator routines of \PETSc{} provide the resulting gradients, 
using the adjoint $v_0$ for the objective with respect to initial conditions
as provided by Equation~\ref{eq:gradient}. 

To compute the gradients, the user first creates the \lstinline{TS} object for a regular forward integration, instructs it to save the trajectories needed for the backward integration,
and then calls \lstinline{TSSolve()}.
The user must provide the initial conditions for the adjoint integration, as computed from Equation~\ref{eq:initadj}, and initiate the reverse mode computation via the following.
\begin{lstlisting}
TSAdjointSolve(TS ts);
\end{lstlisting}

The trajectory information is managed by an extensible abstract trajectory object, \lstinline{TSTrajectory}, that has multiple implementations.
This object provides support for storing checkpoints at particular time steps and recomputing
the missing information. The \lstinline{revolve} \citep{griewank2000algorithm} library is used by
\lstinline{TSTrajectory} to generate an optimal checkpointing schedule that
minimizes the recomputations given a limited number of available
checkpoints; see Section~\ref{sec:checkpointing}. 
\PETSc{} and the \trl{revolve} library also provide an optimal multistage
checkpointing scheme that uses both RAM and disk for storage.
We remark that for higher-order in time integrators the algorithm requires several forward states (or stage values
in the context of multistage time-stepping methods) to evaluate the
Jacobian matrices during the adjoint integration.

\subsection{TAO: gradient-based optimization}
\label{sec:opt}

The Toolkit for Advanced Optimization is a scalable software
library for solving large-scale optimization applications on
high-performance architectures.  It is motivated by the limited
support for scalable parallel computations and the lack of reuse of linear
algebra software in currently available optimization software.  
\TAO{} contains unconstrained minimization, bound-constrained minimization, 
nonlinear complementary, and nonlinear least-squares solvers.
The structure of these problems can differ significantly, but \TAO{} has a 
similar interface to all its solvers.

\TAO{} applications follow a set of procedures for 
solving an optimization problem using gradient-based methods in a way similar to the TS integrators.

\begin{lstlisting}
TaoCreate(MPI_Comm comm, Tao *tao); 
TaoSetInitialVector(Tao tao, Vec x);
TaoSetObjectiveAndGradientRoutine(Tao tao, (*fg)(Tao,Vec,PetscReal*,Vec,void*), void *ctx);
TaoSolve(Tao tao);
\end{lstlisting}

Note that one may select the solver,  as one can set virtually all \PETSc{} and \TAO{} options, 
at runtime by using an options database.  Through this
database, the user can select not only a minimization method but also the convergence
tolerance, various monitoring routines, iterative methods
and preconditioners for solving the linear systems, and so forth.

\TAO{} provides two gradient-based optimization algorithms appropriate for the problems considered in this paper.
The {\em limited-memory, variable-metric} (LMVM) \citep{liu1989limited} method computes a positive-definite
approximation to the Hessian matrix from a limited number of previous
iterates and gradient evaluations.  A direction is obtained by
solving the system of equations
\[
H_k d_k = -\nabla f(x_k),
\]
where $H_k$ is the Hessian approximation obtained by using the Broyden update \citep{NW99}
formula.  The inverse of $H_k$ can readily be applied to obtain the 
direction $d_k$.  Having obtained the direction, a Mor\'{e}-Thuente 
line search \citep{more1994line} is applied to compute a step length, $\tau_k$, that 
approximately solves the one-dimensional optimization problem
\[
\min_\tau f(x_k + \tau d_k).
\]
The current iterate and Hessian approximation are updated, and the process
is repeated until the method converges.  All the numerical studies in this paper utilize LMVM.
Note the approximate Hessian discussed above is computed by the TAO algorithm; the user does not provided it.

For this work, we utilize both the ODE/DAE integrators of TS and the gradient-based
optimization of \TAO{} for PDE-constrained optimization.  For spatial discretization, the user (or an appropriate library) needs to
provide the function evaluations and the associated
Jacobians that arise from the space-dependent operators equipped with
the appropriate boundary conditions. In addition, the user must define
the objective function; however, the rest of the process is handled by the
libraries. 

\section{Computational Results}
The matrix-free discretization using spectral elements is illustrated on the viscous Burgers equation. The time integration and the optimization algorithm are performed using the PETSc capabilities supported by  TS and TAO.
To start, we validate the results against a classical solution and perform a grid convergence study in one dimension. We then present a three-dimensional problem and perform scalability studies, to illustrate both the superiority of the spatial discretization and the efficiency and flexibility of the \PETSc{} software stack.
Since computing with Jacobians is fraught with human error, both in computing the terms mathematically and correctly coding them,
\PETSc{} provides a multitude of ways for checking these terms. We utilized two of these approaches as consistency checks on the algorithms and code in this work. First we
compared the gradients computed by the adjoint approach with finite differences applied to the ODE solver, and second, we computed the Jacobian
explicitly by multiplying by the
identity matrix using the matrix-free Jacobian application and confirming that the transpose product matches those from the explicit matrix representation.
As expected, the results
matched to roughly seven digits indicating the likelihood that the algorithm and code are correct. 

\subsection{Convergence study}
To study the convergence of the implementation, we freeze the nonlinear term of Burgers equation and study the convergence of the advection-diffusion equation.
We utilize the exact reference solution of the form
\begin{equation*}
u(x,t) = \sum_{j=1}^{nc} a_j \sin(2 \pi j(x - a t))e^{-\nu 4 \pi^2 j^2 t},
\end{equation*}
where $ nc $ is 5 and the $ a_j $ are selected from a uniform
distribution between 0.9 and 1. The initial guess is of the same form
with $ t $ equal to zero and utilizing different random values. The
diffusion coefficient is $ \nu = 10^{-5}$, the advection speed is $a=
0.1$, and the short time horizon $T=0.01$. This choice of parameters renders the PDE convection dominated,
which is particularly useful in isolating the ill-posedness stemming from the diffusive component.

 The tolerances for the \TAO{} solver were
set low so that they do not introduce any errors. The time step
was selected adaptively so that the error due to the time
discretization is also small.

\begin{figure}[!htb]
\centering \subfloat[h-refinement]{\includegraphics[width=0.47\textwidth]{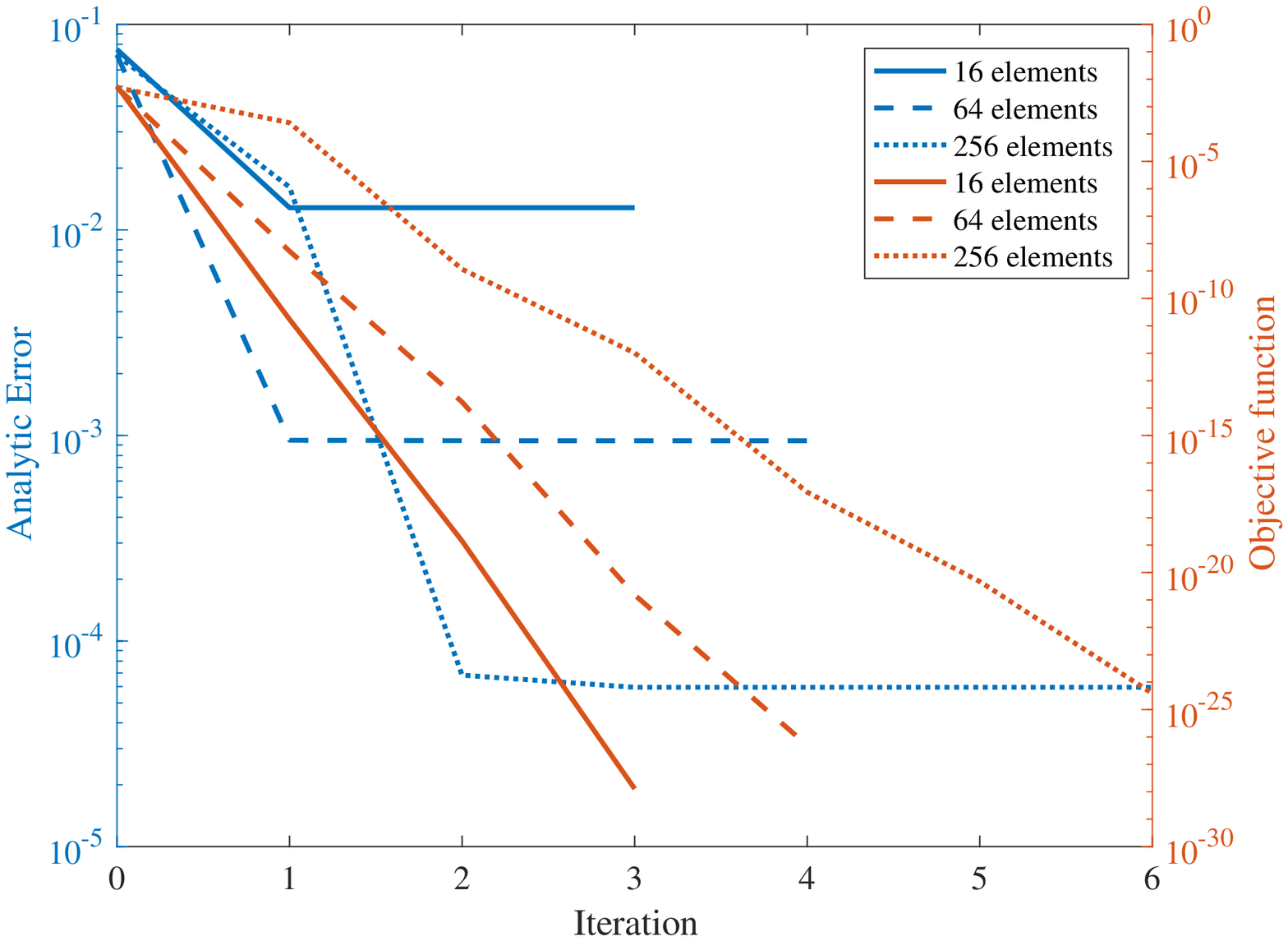}
\label{fig:h-refine}}
\quad \subfloat[p-refinement]{\includegraphics[width=0.47\textwidth]{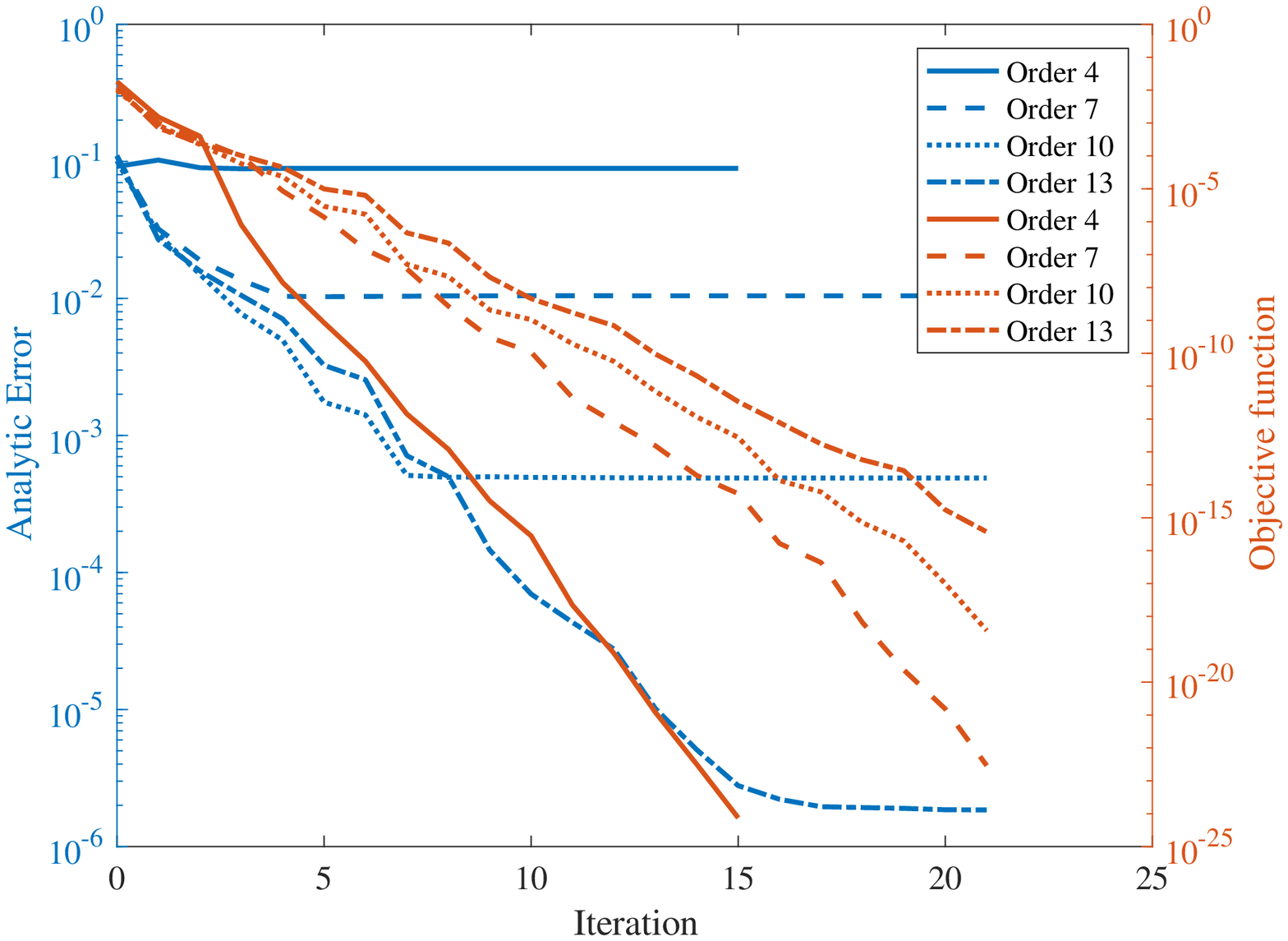}
\label{fig:p-refine}}
\caption{Convergence history of  analytic solution and objective function.}
\label{fig:refine}
\end{figure}

\begin{figure}[!htb]
\centering \subfloat[h-refinement]{\includegraphics[width=0.47\textwidth]{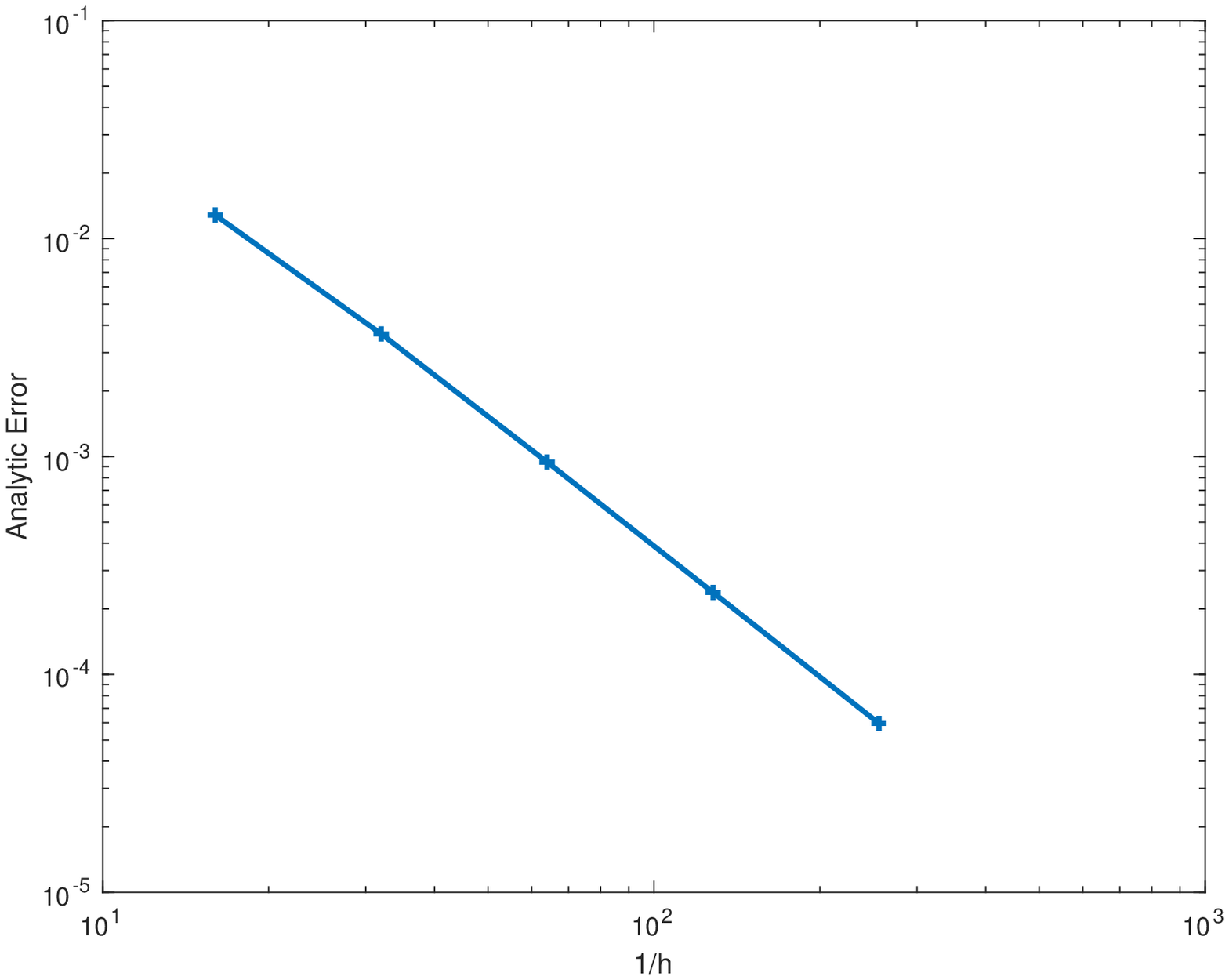}
\label{fig:h-refine2}}
\quad \subfloat[p-refinement]{\includegraphics[width=0.47\textwidth]{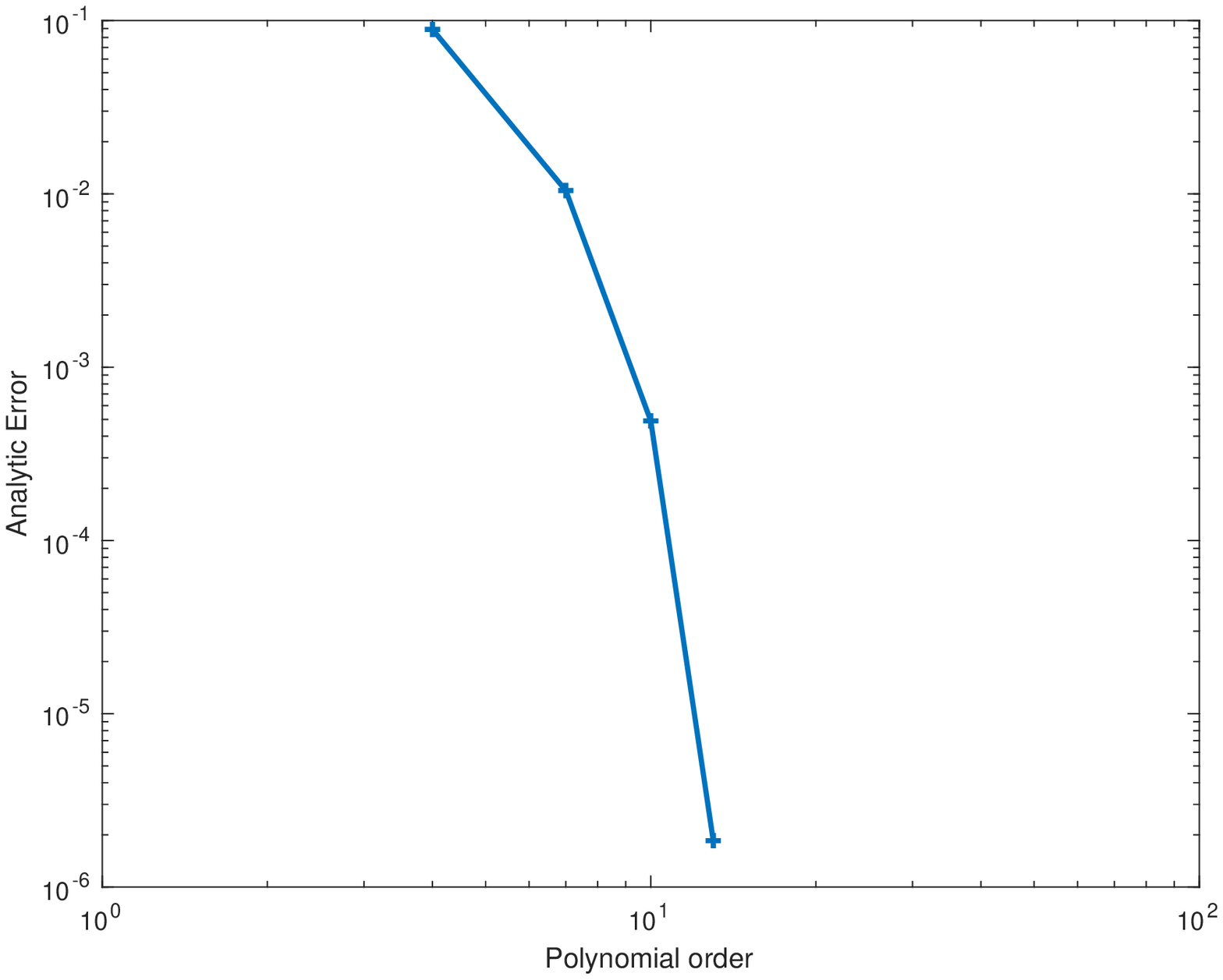}
\label{fig:p-refine2}}
\caption{Analytic error with spatial refinement.}
\label{fig:refine-2}
\end{figure}

In Figure \ref{fig:refine}, we plot the convergence history for the
iterations of the optimization algorithm for a variety of h- and
p-refinements. The $L^2$ norm of the error to the analytic optimization problem is
plotted in Figure \ref{fig:refine-2}. As expected, the convergence is
second-order for h-refinement and spectral for p-refinement.

\subsubsection{Viscous Burgers equation}

The time-dependent viscous Burgers equation in one dimension is given 
by the partial differential equation
\[
u_t -  \nu \Delta u + u \nabla u = 0.
\]

We follow the work of \citep{ou_2011},  which was also performed in the 
framework of high-order spectral discretizations, albeit based on a Chebyshev grid and using 
discontinuous Galerkin; however, the authors restrict themselves to low-order polynomials, and their 
optimization algorithm does not exhibit the same acceleration that we achieve.

The idea is once again to choose an analytical solution and seek to converge to it in $L^2$ norm.
The analytical solution is
\begin{equation}
u(x,t)=\frac{2\nu\pi\sin(\pi x)e^{-\nu t \pi^2}}{2+e^{-\nu t \pi^2}\cos(\pi x)},
\end{equation}
where the exact initial condition can be obtained for $t=0$ and the time horizon set for $t=T$.
We start from an initial condition that deviates from the analytical one by 
a Gaussian signal centered in the middle of the domain of interest $[-2,\ 2]$:
\begin{equation}
u(x,0)=\frac{2\nu\pi\sin(\pi x)}{2+\cos(\pi x)}+e^{-4(x-2)^2}.
\end{equation}

\begin{figure}[!htb]
\centering
\subfloat[Burgers equation, T=4 time horizon, optimal solution (red),  initial condition at first iteration (black)]
{\includegraphics[width=0.47\textwidth]{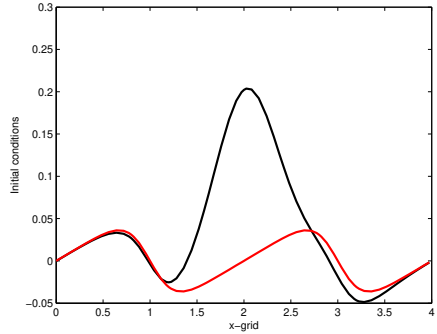}
\label{fig:ic_obj_burgers}}
\quad
\subfloat[Objective function and initial condition error decay versus number of iterations;  Data from \citep{ou_2011} (red), 
 error between computed initial condition and analytic initial condition (blue),  objective function error (black).]
{\includegraphics[width=0.47\textwidth]{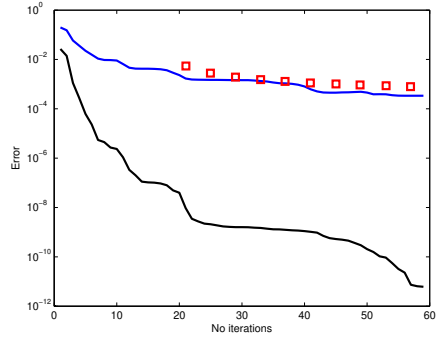}
\label{fig:decay_burgers}}
\caption{Optimal initial conditions for Burgers equation with an analytic
 reference function up to time horizon T=4.}
\end{figure}

\begin{figure}[!htb]
\centering
\subfloat[Initial condition and evolution towards= the analytic initial condition.]
{\includegraphics[width=0.47\textwidth]{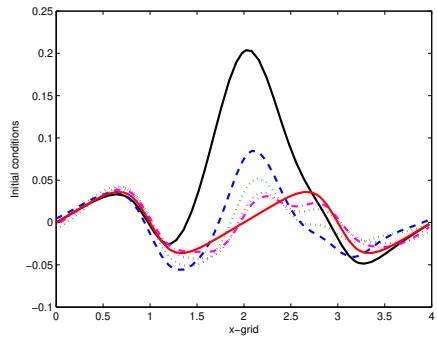}
\label{fig:ic_burgers}}
\quad
\subfloat[Solution of Burgers equation and evolution toward the objective function.]
{\includegraphics[width=0.47\textwidth]{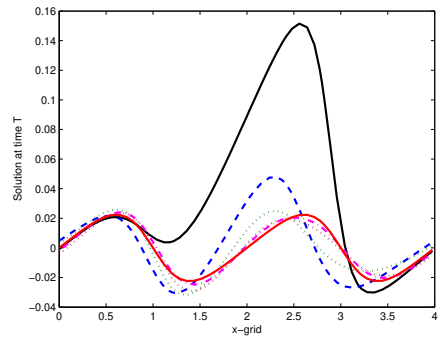}
\label{fig:final_burgers}}
\caption{Convergence of the solution for the first 4 iterations: (black) iteration 0, (blue) iteration 1, (green) iteration 2, (brown) iteration 3, (purple) iteration 4.}
\end{figure}

\begin{figure}[!htb]
\centering
{\includegraphics[width=0.47\textwidth]{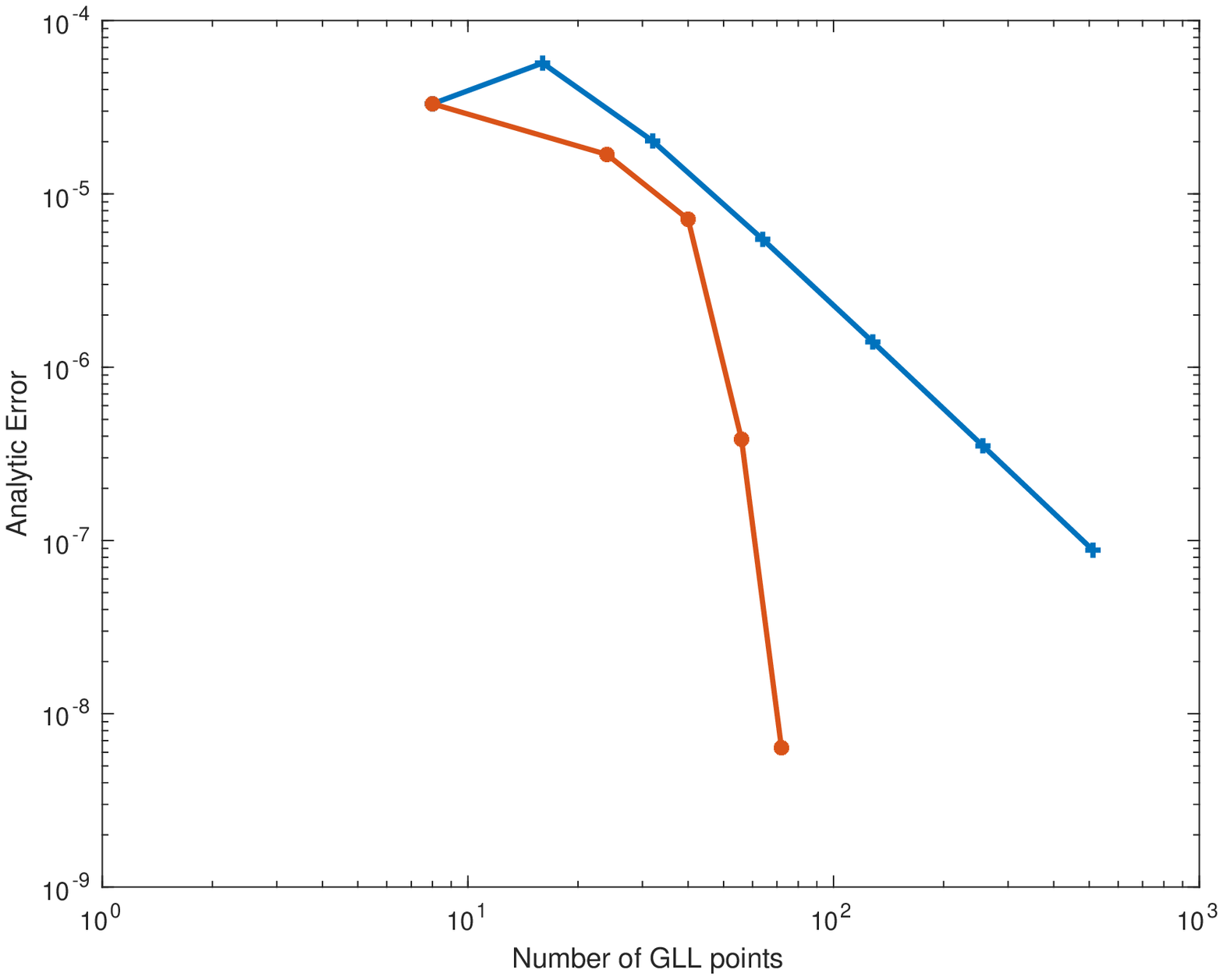}}
\caption{Convergence to the analytic optimal solution with h-refinement (blue) and 
p-refinement (red) for Burgers equation.\label{fig:bconv}.}
\end{figure}
In Figure~\ref{fig:bconv} we study the h- and p-refinement convergence for Burgers equation with a $\nu$ of 0.001. As expected, there is quadratic convergence with h-refinement and spectral with p-refinement. For this
case we had to use the analytic initial conditions to obtain convergence to the discrete solution; 
using the previous perturbed Gaussian as initial conditions results in stagnation.

 \subsection{Three-dimensional viscous Burgers equation}
 For modeling simplicity
 we compute with periodic boundary conditions while using initial conditions given by
 \begin{eqnarray}
 u_1(\mathbf x,0)&=& \sin(0.5\pi x_1)\cdot\cos(0.5\pi x_2)\cdot\cos(0.5\pi x_3)\\ \nonumber
 u_2(\mathbf x,0)&=& \sin(0.5\pi x_2)\cdot\cos(0.5\pi x_3)\cdot\cos(0.5\pi x_1)\\ \nonumber
 u_3(\mathbf x,0)&=& \sin(0.5\pi x_3)\cdot\cos(0.5\pi
 x_1)\cdot\cos(0.5\pi x_2) \,.
 \end{eqnarray}
 The objective function is the initial condition pre-multiplied by $\mathrm e^{-\nu\cdot t}$. This would provide an analytical solution if only diffusion were present in the process. Although this is not the case, it provides a suitable solution to the optimization problem. We note that the objective function needs to be a solution to the partial differential equation. Additionally, we seek short time horizons that do not allow the development of shocks, unavoidable due to the nonlinearity of the advection term.
 
 \begin{figure}[!htb]
\centering
\subfloat[Desired solution, 3D slice through $u_2$]
{\includegraphics[width=0.47\textwidth]{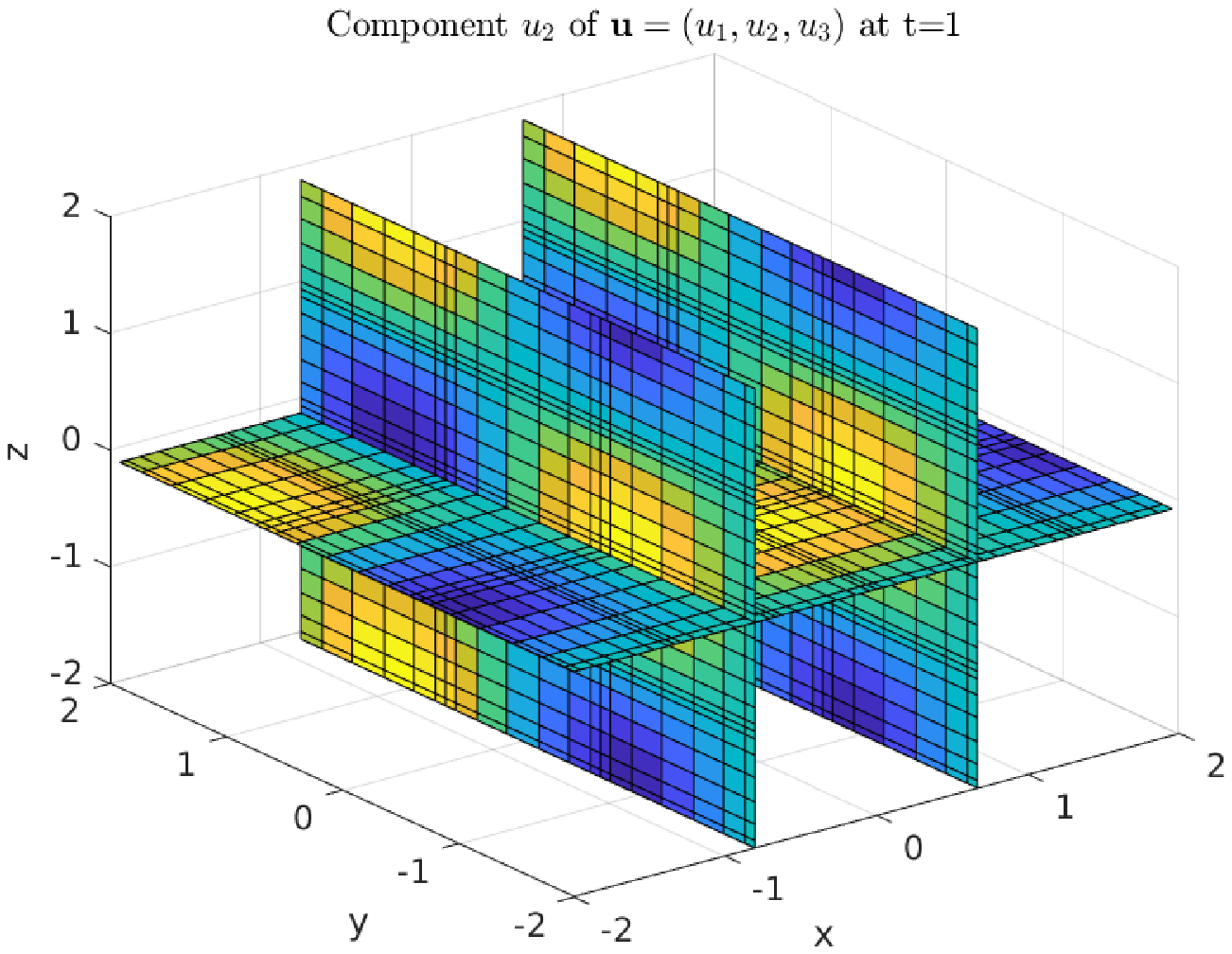}
\label{fig:obj3d}}
\quad
\subfloat[Initial condition $u_2$ through midplane $z=0$]
{\includegraphics[width=0.47\textwidth]{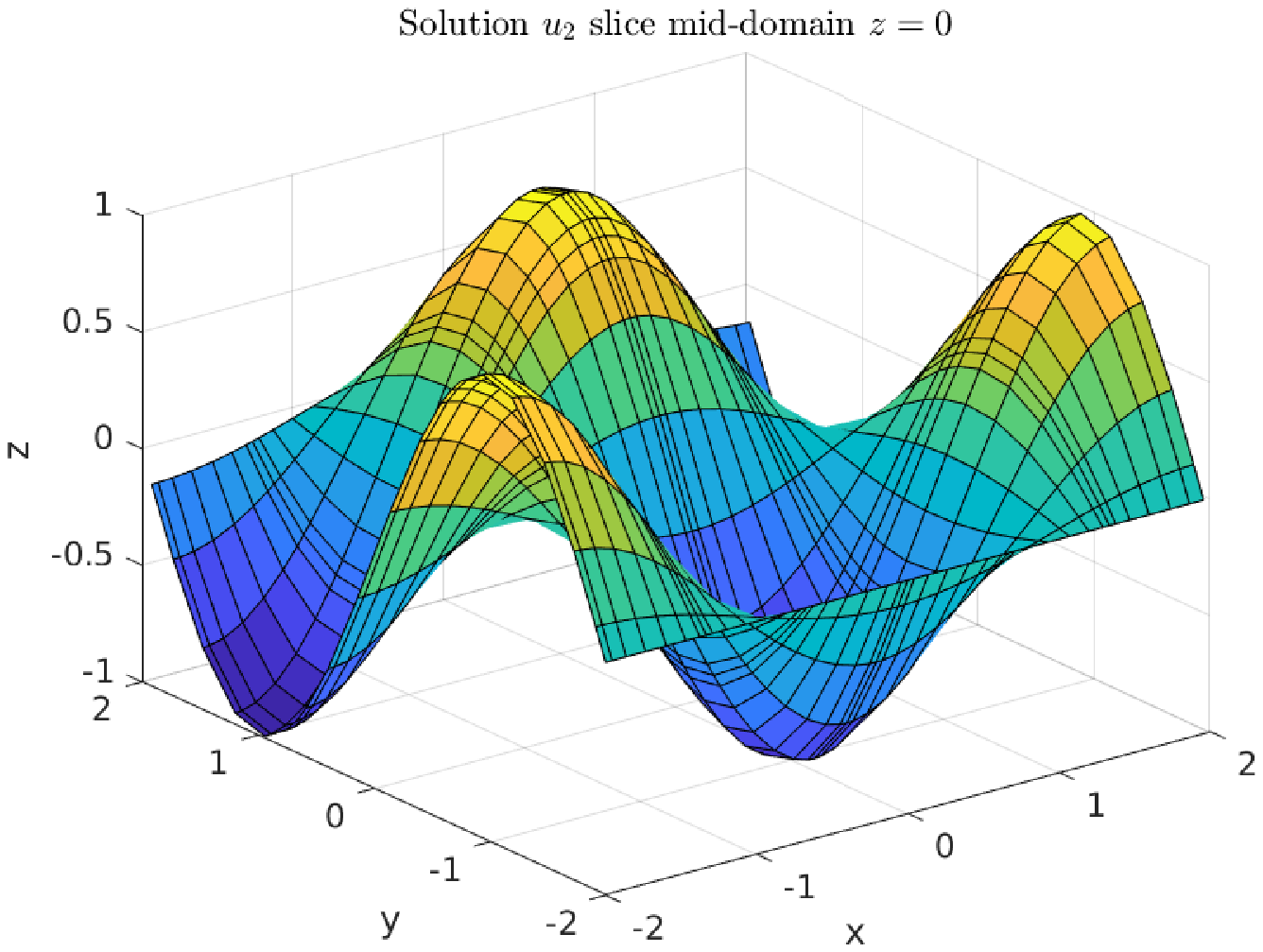}
\label{fig:init3d}}
\caption{Three-dimensional Burgers.}
\end{figure}

The objective function and initial conditions are illustrated in Figure~\ref{fig:obj3d}, and Figure~\ref{fig:init3d}, respectively, for the field $u_2$.

\subsection{Scalability studies}
 The scalability studies were performed on the Argonne Leadership Computing Facility Theta computer, which consists of Intel KNL (second generation) nodes, each with 64 compute cores. Intel's MKL BLAS was used for the element matrix-matrix products. We consider two cases: a small case with $16^3$
 elements of identical degree $P=8$, resulting in a problem with 4,214,784 degrees of freedom, and a larger case of $ 64^3$ elements of identical degree resulting in 269,746,176 unknowns.
 We use two integration schemes: the explicit Runge-Kutta (RK-3) with three stages and the implicit Crank-Nicolson (CN) scheme. For timing,
 the runs were made for five optimization steps, which resulted in a decrease in the initial objective function value of around $2 \times 10^{-3}$.

\begin{figure}[!htb]
\centering
\subfloat[RK-3 integrator (explicit)]
{\includegraphics[trim = 0.1in 2.5in 0.1in 2.5in, clip, width=0.47\textwidth]{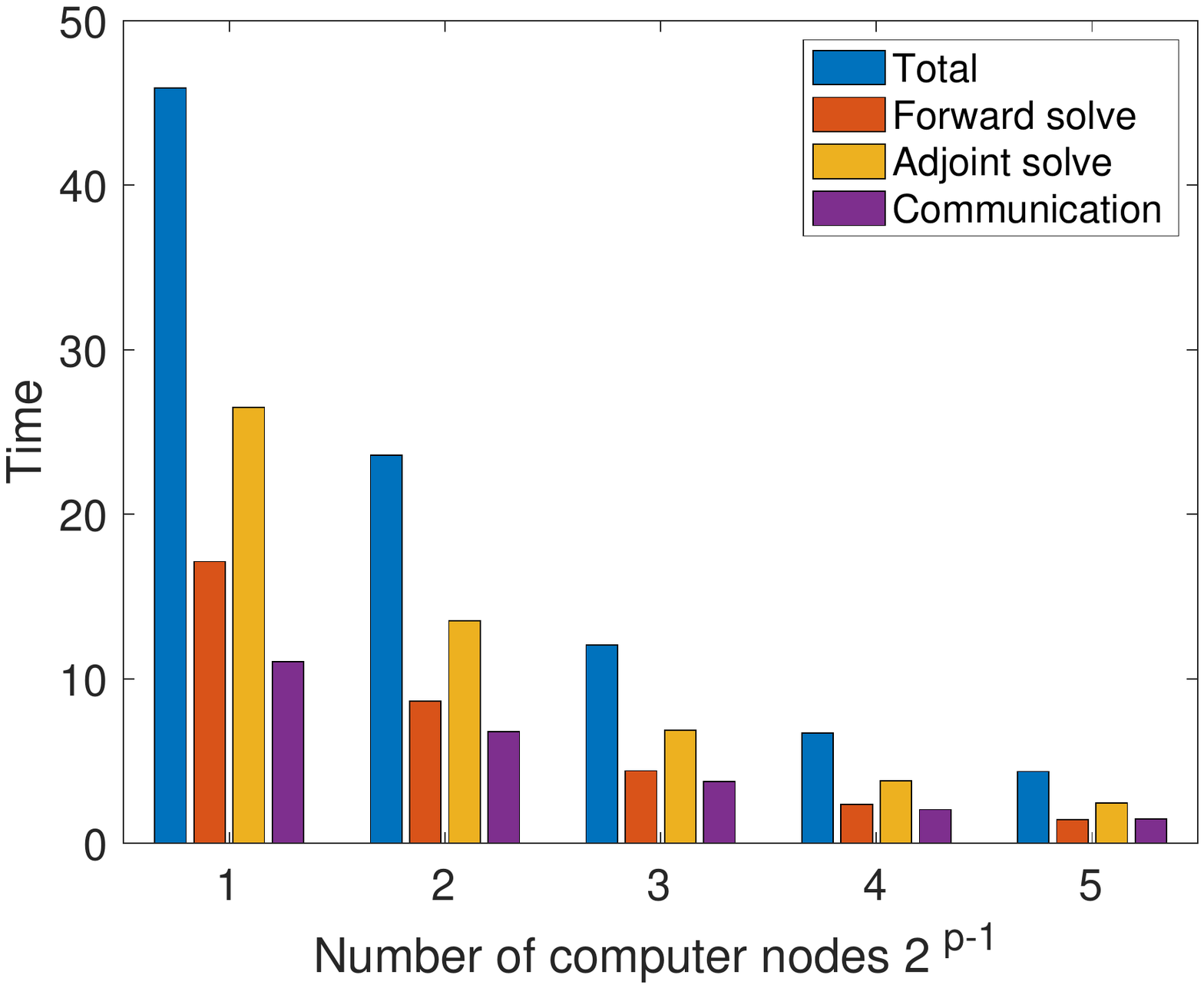}
\label{fig:bar_rk}}
\quad
\subfloat[Crank-Nicolson integrator (implicit)]
{\includegraphics[trim = 0.1in 2.5in 0.1in 2.5in, clip, width=0.47\textwidth]{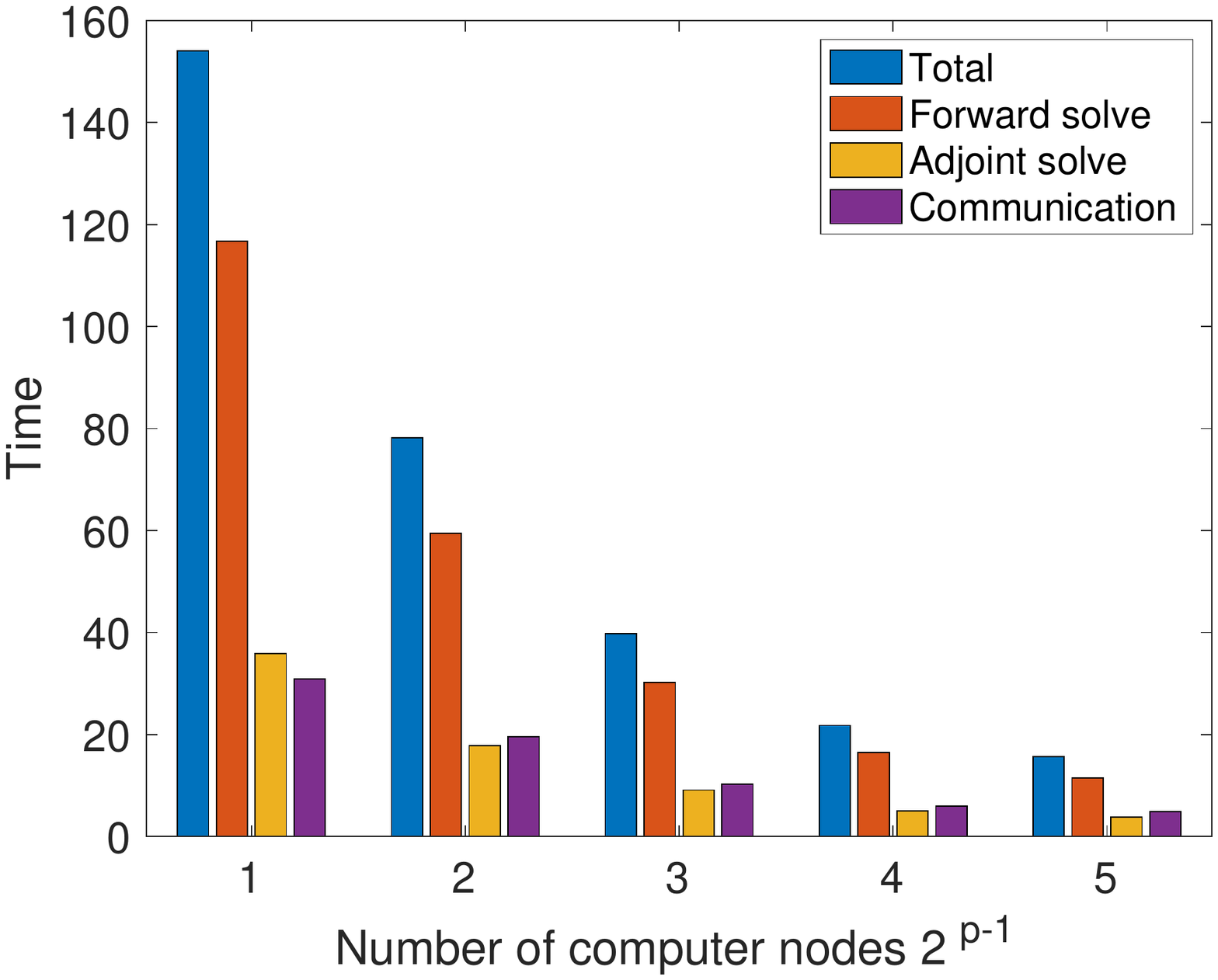}
\label{fig:bar_cn}}
\caption{Comparison of time needed for explicit and implicit forward and adjoint integrators for the small problem.}
\end{figure}

 In the first study, we do not use revolve but instead store all the history in the large ``slow'' memory available on KNL systems. 
 In Figures \ref{fig:bar_rk} and \ref{fig:bar_cn}, we show the total time spent in the solver as well as the time spent in the forward integration and the adjoint integration for the small problem. In addition, we show the time spent in the communication (which occurs within both of the time integrators). For the explicit integrator, as expected,
 slightly more time is needed for the adjoint integration than for the forward integration. With the implicit integrator, the forward simulation requires significantly more time than does the adjoint integration since the adjoint integrator requires only linear solves whereas the forward integrations require nonlinear solvers.

\begin{figure}[!htb]
\centering
\subfloat[RK-3 integrator (explicit)]
{\includegraphics[trim = 0.1in 2.5in 0.1in 2.5in, clip, width=0.47\textwidth]{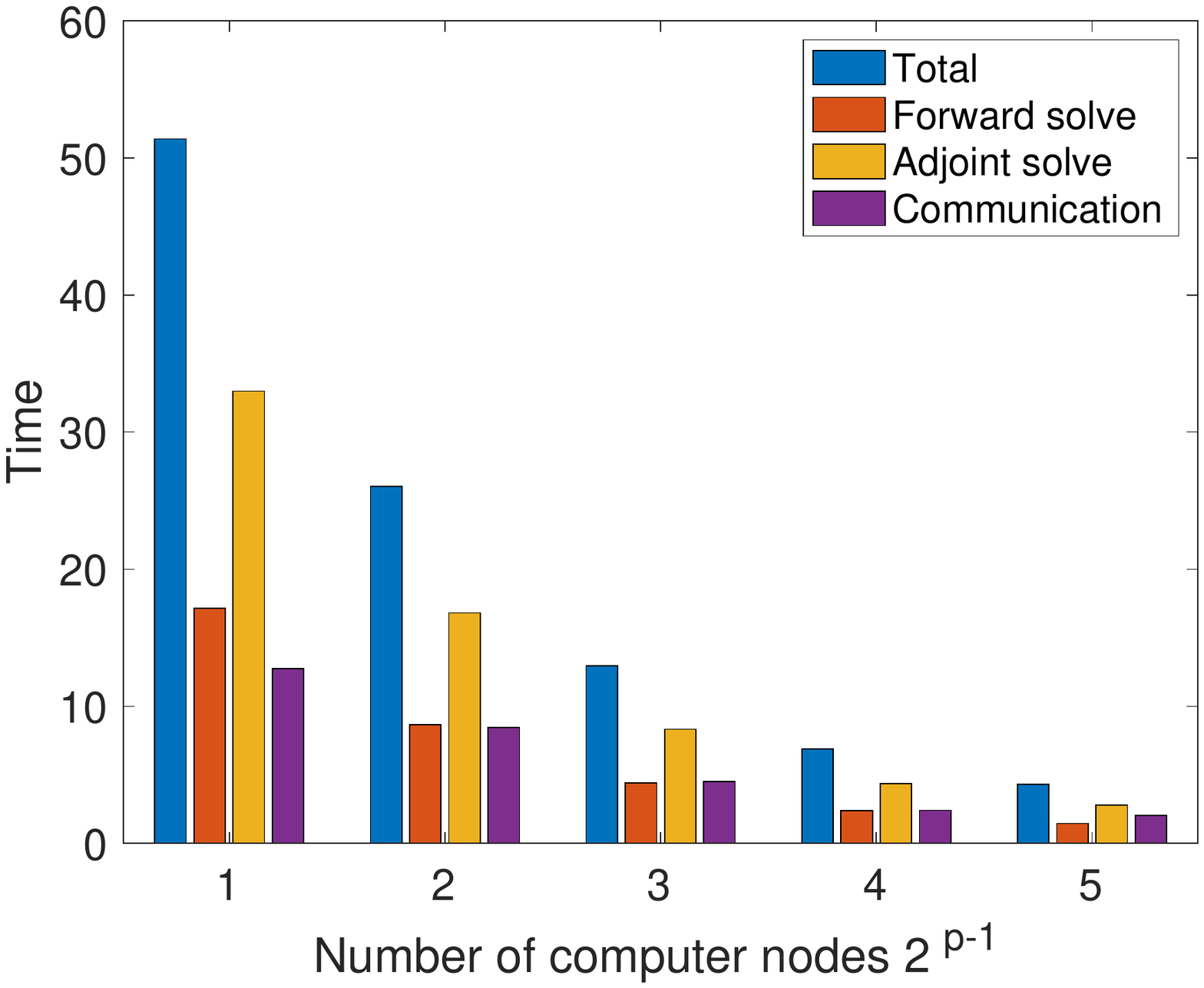}
\label{fig:bar_rk_revolve}}
\quad
\subfloat[Crank-Nicolson integrator (implicit)]
{\includegraphics[trim = 0.1in 2.5in 0.1in 2.5in, clip, width=0.47\textwidth]{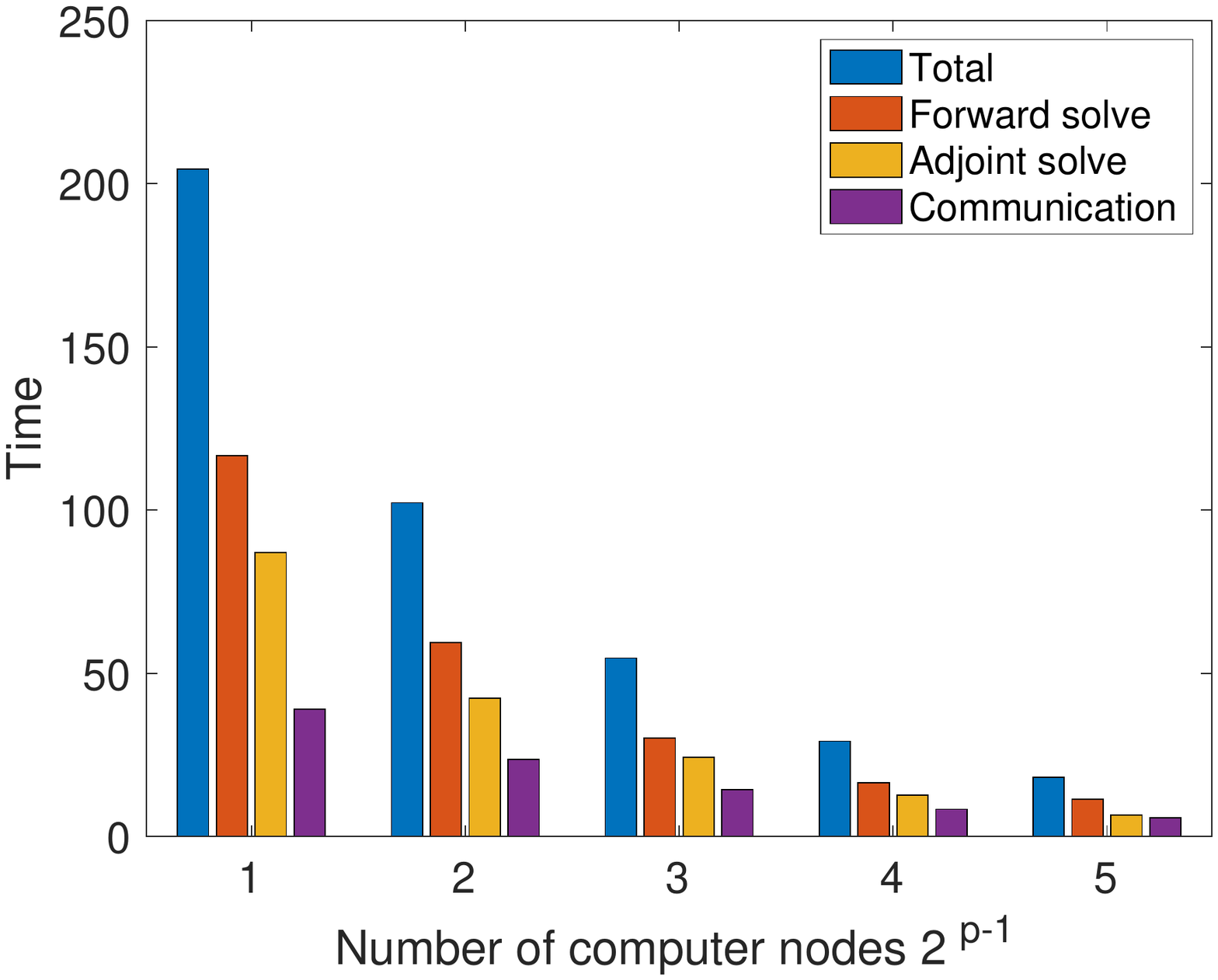}
\label{fig:bar_cn_revolve}}
\caption{Comparison of time needed for explicit and implicit forward and adjoint integrators for the small problem using revolve with a storage size of 10 vectors.}
\end{figure}

 In Figures \ref{fig:bar_rk_revolve} and \ref{fig:bar_cn_revolve} we show the time needed for the forward and adjoint integration when revolve is used with a storage space
 size of 10 vectors.
 Note that the additional time needed to recompute portions of the forward integration is attributed to the adjoint integration.
 
\begin{figure}[!htb]
\centering
\includegraphics[trim = 0.1in 2.5in 0.1in 2.5in, clip, width=0.5\textwidth]{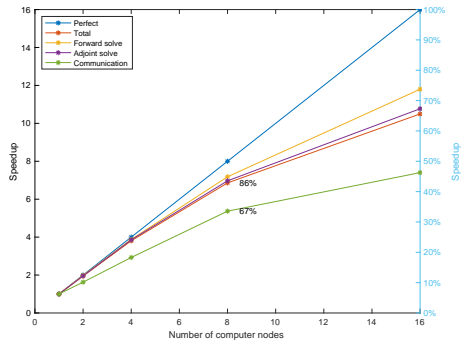}
\caption{Speedup for small problem RK-3.}\label{fig:rk_speedup}
\end{figure}

Figure \ref{fig:rk_speedup} shows the speedup for the small problem solved with RK-3. Note that the communication cost dominates the poor performance for larger process counts. 
For 16 compute nodes, there are 1,024 cores
and only 4,116 unknowns per process; this limits the speedup. The application of the nonlinear function needed by TS for the forward integration
requires one global to local communication (providing the ghost points needed for the element computations) and one local to global communication to accumulate the
results back into the global vector. The application of the transpose Jacobian product requires three communications, one global to local for the vector at which
the Jacobian is computed, one for the transpose matrix-vector product input value, and one to accumulate the result.
Thus, even though the adjoint application is linear while the forward application is nonlinear, the adjoint application requires more communication.

\begin{figure}[!htb]
\centering
\subfloat[RK-3 integrator (explicit)]
{\includegraphics[trim = 0.1in 2.5in 0.1in 2.5in, clip, width=0.47\textwidth]{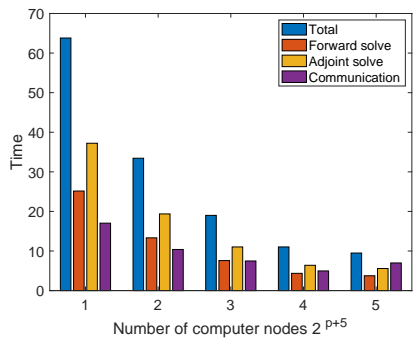}
\label{fig:bar_rk_large}}
\quad
\subfloat[RK-3 integrator (explicit)]
{\includegraphics[trim = 0.1in 2.5in 0.1in 2.5in, clip, width=0.47\textwidth]{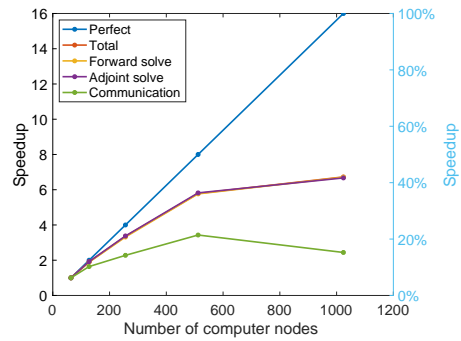}
\label{fig:rk_speedup_large}}
\caption{Comparison of time needed and speedup for the explicit forward and adjoint integration for the large problem.}
\end{figure}

Figures \ref{fig:bar_rk_large} and \ref{fig:rk_speedup_large} repeat
the results of the computations for the large problem using RK-3. It
is run from 64 nodes (4,096 cores) to 1,024 nodes (65,536 cores) on a
problem with 269,746,176 unknowns. Note that for the largest number of
computer nodes, we again have the extremely modest 4,116 degrees of freedom per compute
process, which amounts to four elements.  Again, the communication time dominates the solution
time; the ALCF Intel KNL system is known for having slow communication
time relative to floating-point performance.

 Our final results are for 2,048 nodes (131,072 cores) and 2,157,969,408 unknowns shown in Table \ref{table:extralarge}. Again the adjoint integration requires
 a bit more time than the forward integration and the communication time dominates. 
\begin{table}[!htb]
\begin{tabular}{c c c c}
\hline
Optimization & Forward integration & Adjoint integration & Communication \\
\hline
20.46 & 8.27 & 11.81 & 7.95 \\
\hline
\end{tabular}
\caption{Timing (seconds) results using RK-3 for the extra-large problem.}
\label{table:extralarge}
\end{table}

\section{Conclusions}
\label{sec:conclusion}
We have demonstrated how \PETSc{} time integrators, it's adjoint
capability, and the gradient-based optimization capabilities of \TAO{}
can be used to solve PDE-constrained optimization systematically
problems using the spectral element method with both explicit and implicit time integrators. We provided a new
efficient approach to apply the transpose of the Jacobian using tensor
contractions for the nonlinear Burgers equation, thus resulting in an adjoint integration time
that is not much more expensive than the forward integration (for explicit integrators). In addition,
we demonstrated the performance and scalability of the algorithm and code on problems up to
with over 2 billion unknowns on over 130,000 compute cores. Examples that demonstrate the TSAdjoint and \TAO{} capabilities of \PETSc{} are available in the \PETSc{}
git repository at \href{https://gitlab.com/petsc/petsc}{https://gitlab.com/petsc/petsc} in the directory {\em src/ts/examples/tutorials} and its subdirectories.

\section*{Acknowledgments}

We thank Hong Zhang for the time-integration adjoint capability implemented in \PETSc{}.  This material is based upon work supported by the
U.S. Department of Energy, Office of Science, Advanced Scientific
Computing Research under Contract DE-AC02-06CH11357 and the Exascale
Computing Project (Contract No. 17-SC-20-SC).  This research was
supported by the Exascale Computing Project (17-SC-20-SC), a
collaborative effort of the U.S. Department of Energy Office of
Science and the National Nuclear Security Administration.

\bibliographystyle{elsarticle-harv}

\begin{thebibliography}{34}
\expandafter\ifx\csname natexlab\endcsname\relax\def\natexlab#1{#1}\fi
\providecommand{\url}[1]{\texttt{#1}}
\providecommand{\href}[2]{#2}
\providecommand{\path}[1]{#1}
\providecommand{\DOIprefix}{doi:}
\providecommand{\ArXivprefix}{arXiv:}
\providecommand{\URLprefix}{URL: }
\providecommand{\Pubmedprefix}{pmid:}
\providecommand{\doi}[1]{\href{http://dx.doi.org/#1}{\path{#1}}}
\providecommand{\Pubmed}[1]{\href{pmid:#1}{\path{#1}}}
\providecommand{\bibinfo}[2]{#2}
\ifx\xfnm\relax \def\xfnm[#1]{\unskip,\space#1}\fi
\bibitem[{Abhyankar et~al.(2014)Abhyankar, Brown, Constantinescu, Ghosh and
  Smith}]{tspaper}
\bibinfo{author}{Abhyankar, S.}, \bibinfo{author}{Brown, J.},
  \bibinfo{author}{Constantinescu, E.}, \bibinfo{author}{Ghosh, D.},
  \bibinfo{author}{Smith, B.F.}, \bibinfo{year}{2014}.
\newblock \bibinfo{title}{{PETSc/TS}: A Modern Scalable {DAE/ODE} Solver
  Library}.
\newblock \bibinfo{type}{Preprint} \bibinfo{number}{ANL/MCS-P5061-0114}. ANL.
\bibitem[{Amir and Sigmund(2011)}]{amir2011reducing}
\bibinfo{author}{Amir, O.}, \bibinfo{author}{Sigmund, O.},
  \bibinfo{year}{2011}.
\newblock \bibinfo{title}{On reducing computational effort in topology
  optimization: {H}ow far can we go?}
\newblock \bibinfo{journal}{Structural and Multidisciplinary Optimization}
  \bibinfo{volume}{44}, \bibinfo{pages}{25--29}.
\bibitem[{Balay et~al.(2020)Balay, Abhyankar, Adams, Brown, Brune, Buschelman,
  Dalcin, Dener, Eijkhout, Gropp, Karpeyev, Kaushik, Knepley, May, McInnes,
  Mills, Munson, Rupp, Sanan, Smith, Zampini, Zhang and Zhang}]{petsc-user-ref}
\bibinfo{author}{Balay, S.}, \bibinfo{author}{Abhyankar, S.},
  \bibinfo{author}{Adams, M.F.}, \bibinfo{author}{Brown, J.},
  \bibinfo{author}{Brune, P.}, \bibinfo{author}{Buschelman, K.},
  \bibinfo{author}{Dalcin, L.}, \bibinfo{author}{Dener, A.},
  \bibinfo{author}{Eijkhout, V.}, \bibinfo{author}{Gropp, W.D.},
  \bibinfo{author}{Karpeyev, D.}, \bibinfo{author}{Kaushik, D.},
  \bibinfo{author}{Knepley, M.G.}, \bibinfo{author}{May, D.A.},
  \bibinfo{author}{McInnes, L.C.}, \bibinfo{author}{Mills, R.T.},
  \bibinfo{author}{Munson, T.}, \bibinfo{author}{Rupp, K.},
  \bibinfo{author}{Sanan, P.}, \bibinfo{author}{Smith, B.F.},
  \bibinfo{author}{Zampini, S.}, \bibinfo{author}{Zhang, H.},
  \bibinfo{author}{Zhang, H.}, \bibinfo{year}{2020}.
\newblock \bibinfo{title}{{PETS}c Users Manual}.
\newblock \bibinfo{type}{Technical Report} \bibinfo{number}{ANL-95/11 -
  Revision 3.13}. Argonne National Laboratory.
\bibitem[{Brown et~al.(2019)Brown, Abdelfattah, Barra, Dobrev, Dudouit,
  Fischer, Kolev, Medina, Min, Ratnayaka, Smith, Thompson, Tomov, Tomov and
  Warburton}]{ceed_report}
\bibinfo{author}{Brown, J.}, \bibinfo{author}{Abdelfattah, A.},
  \bibinfo{author}{Barra, V.}, \bibinfo{author}{Dobrev, V.},
  \bibinfo{author}{Dudouit, Y.}, \bibinfo{author}{Fischer, P.},
  \bibinfo{author}{Kolev, T.}, \bibinfo{author}{Medina, D.},
  \bibinfo{author}{Min, M.}, \bibinfo{author}{Ratnayaka, T.},
  \bibinfo{author}{Smith, C.}, \bibinfo{author}{Thompson, J.},
  \bibinfo{author}{Tomov, S.}, \bibinfo{author}{Tomov, V.},
  \bibinfo{author}{Warburton, T.}, \bibinfo{year}{2019}.
\newblock \bibinfo{title}{{CEED ECP Milestone Report: Public release of CEED
  2.0}}.
\newblock \URLprefix \url{https://doi.org/10.5281/zenodo.2641316},
  \DOIprefix\doi{10.5281/zenodo.2641316}.
\bibitem[{Constantinescu(2018)}]{constantinescu2015estimating}
\bibinfo{author}{Constantinescu, E.}, \bibinfo{year}{2018}.
\newblock \bibinfo{title}{Generalizing global error estimation for ordinary
  differential equations by using coupled time-stepping methods}.
\newblock \bibinfo{journal}{Journal of Computational and Applied Mathematics}
  \bibinfo{volume}{332}, \bibinfo{pages}{140--158}.
\newblock \URLprefix
  \url{http://www.sciencedirect.com/science/article/pii/S0377042717302480},
  \DOIprefix\doi{https://doi.org/10.1016/j.cam.2017.05.012}.
\bibitem[{Dener et~al.(2019)Dener, Denchfield, Munson, Sarich, Wild, Benson and
  McInnes}]{tao-user-ref}
\bibinfo{author}{Dener, A.}, \bibinfo{author}{Denchfield, A.},
  \bibinfo{author}{Munson, T.}, \bibinfo{author}{Sarich, J.},
  \bibinfo{author}{Wild, S.}, \bibinfo{author}{Benson, S.},
  \bibinfo{author}{McInnes, L.C.}, \bibinfo{year}{2019}.
\newblock \bibinfo{title}{{Toolkit for Advanced Optimization (TAO) users
  manual}}.
\newblock \bibinfo{type}{Technical Report} \bibinfo{number}{ANL/MCS-TM-322 -
  Rev 3.12}. Argonne National Laboratory.
\bibitem[{Deville et~al.(2002)Deville, Fischer and Mund}]{fischer:hom}
\bibinfo{author}{Deville, M.O.}, \bibinfo{author}{Fischer, P.F.},
  \bibinfo{author}{Mund, E.H.}, \bibinfo{year}{2002}.
\newblock \bibinfo{title}{High-Order Methods for Incompressible Fluid Flow}.
\newblock \bibinfo{publisher}{Cambridge University Press}.
\newblock \URLprefix \url{http://dx.doi.org/10.1017/CBO9780511546792}.
  \bibinfo{note}{cambridge Books Online}.
\bibitem[{Dunning et~al.(2017)Dunning, Huchette and Lubin}]{dunning2017jump}
\bibinfo{author}{Dunning, I.}, \bibinfo{author}{Huchette, J.},
  \bibinfo{author}{Lubin, M.}, \bibinfo{year}{2017}.
\newblock \bibinfo{title}{{JuMP}: a modeling language for mathematical
  optimization}.
\newblock \bibinfo{journal}{SIAM Review} \bibinfo{volume}{59},
  \bibinfo{pages}{295--320}.
\bibitem[{Falgout(2017)}]{hypre-users-manual}
\bibinfo{author}{Falgout, R.}, \bibinfo{year}{2017}.
\newblock \bibinfo{title}{{hypre} Users Manual}.
\newblock \bibinfo{type}{Technical Report} \bibinfo{number}{Revision 2.11.2}.
  Lawrence Livermore National Laboratory.
\bibitem[{Farrell et~al.(2013)Farrell, Ham, Funke and
  Rognes}]{farrell2013automated}
\bibinfo{author}{Farrell, P.E.}, \bibinfo{author}{Ham, D.A.},
  \bibinfo{author}{Funke, S.W.}, \bibinfo{author}{Rognes, M.E.},
  \bibinfo{year}{2013}.
\newblock \bibinfo{title}{Automated derivation of the adjoint of high-level
  transient finite element programs}.
\newblock \bibinfo{journal}{SIAM Journal on Scientific Computing}
  \bibinfo{volume}{35}, \bibinfo{pages}{C369--C393}.
\bibitem[{Fischer et~al.(2015)Fischer, Lottes, Kerkemeier, Marin, Heisey,
  Obabko, Merzari and Peet}]{nekdoc}
\bibinfo{author}{Fischer, P.}, \bibinfo{author}{Lottes, J.},
  \bibinfo{author}{Kerkemeier, S.}, \bibinfo{author}{Marin, O.},
  \bibinfo{author}{Heisey, K.}, \bibinfo{author}{Obabko, A.},
  \bibinfo{author}{Merzari, E.}, \bibinfo{author}{Peet, Y.},
  \bibinfo{year}{2015}.
\newblock \bibinfo{title}{{Nek5000: User's manual}}.
\newblock \bibinfo{type}{Technical Report} \bibinfo{number}{ANL/MCS-TM-351}.
  Argonne National Laboratory.
\bibitem[{Gander et~al.(2014)Gander, Kwok and Wanner}]{gander2014constrained}
\bibinfo{author}{Gander, M.J.}, \bibinfo{author}{Kwok, F.},
  \bibinfo{author}{Wanner, G.}, \bibinfo{year}{2014}.
\newblock \bibinfo{title}{Constrained optimization: From {L}agrangian mechanics
  to optimal control and {PDE} constraints}, in:
  \bibinfo{booktitle}{Optimization with PDE Constraints}.
  \bibinfo{publisher}{Springer}, pp. \bibinfo{pages}{151--202}.
\bibitem[{Giles and Pierce(1997)}]{giles1997adjoint}
\bibinfo{author}{Giles, M.B.}, \bibinfo{author}{Pierce, N.A.},
  \bibinfo{year}{1997}.
\newblock \bibinfo{title}{Adjoint equations in {CFD}: duality, boundary
  conditions and solution behaviour}.
\newblock \bibinfo{journal}{AIAA paper} \bibinfo{volume}{1850},
  \bibinfo{pages}{1997}.
\bibitem[{Griewank and Walther(2000)}]{griewank2000algorithm}
\bibinfo{author}{Griewank, A.}, \bibinfo{author}{Walther, A.},
  \bibinfo{year}{2000}.
\newblock \bibinfo{title}{Algorithm 799: revolve: an implementation of
  checkpointing for the reverse or adjoint mode of computational
  differentiation}.
\newblock \bibinfo{journal}{ACM Transactions on Mathematical Software}
  \bibinfo{volume}{26}, \bibinfo{pages}{19--45}.
\bibitem[{Gunzburger(2002)}]{gunzburger2002perspectives}
\bibinfo{author}{Gunzburger, M.D.}, \bibinfo{year}{2002}.
\newblock \bibinfo{title}{Perspectives in flow control and optimization}.
\newblock \bibinfo{publisher}{SIAM}.
\bibitem[{Haber and Hanson(2007)}]{Haber_2007}
\bibinfo{author}{Haber, E.}, \bibinfo{author}{Hanson, L.},
  \bibinfo{year}{2007}.
\newblock \bibinfo{title}{Model problems in {PDE}-constrained optimization}.
\newblock \bibinfo{type}{Technical Report}. TS-200X-00X-A, Emory University.
\bibitem[{Hindmarsh et~al.(2005)Hindmarsh, Brown, Grant, Lee, Serban, Shumaker
  and Woodward}]{sundials05}
\bibinfo{author}{Hindmarsh, A.}, \bibinfo{author}{Brown, P.},
  \bibinfo{author}{Grant, K.}, \bibinfo{author}{Lee, S.},
  \bibinfo{author}{Serban, R.}, \bibinfo{author}{Shumaker, D.},
  \bibinfo{author}{Woodward, C.}, \bibinfo{year}{2005}.
\newblock \bibinfo{title}{{SUNDIALS:} suite of nonlinear and
  differential/algebraic equation solvers}.
\newblock \bibinfo{journal}{ACM Transactions on Mathematical Software}
  \bibinfo{volume}{31}, \bibinfo{pages}{363--396}.
\bibitem[{Isik(2013)}]{isik2013spin}
\bibinfo{author}{Isik, O.R.}, \bibinfo{year}{2013}.
\newblock \bibinfo{title}{Spin up problem and accelerating convergence to
  steady state}.
\newblock \bibinfo{journal}{Applied Mathematical Modelling}
  \bibinfo{volume}{37}, \bibinfo{pages}{3242--3253}.
\bibitem[{Kopriva(2009)}]{kopriva2009}
\bibinfo{author}{Kopriva, D.}, \bibinfo{year}{2009}.
\newblock \bibinfo{title}{Implementing spectral methods for partial
  differential equations: algorithms for scientists and engineers}.
\newblock \bibinfo{publisher}{Springer Science \& Business Media}.
\bibitem[{Li and Demmel(2003)}]{lidemmel03}
\bibinfo{author}{Li, X.S.}, \bibinfo{author}{Demmel, J.W.},
  \bibinfo{year}{2003}.
\newblock \bibinfo{title}{{SuperLU\_DIST}: a scalable distributed-memory sparse
  direct solver for unsymmetric linear systems}.
\newblock \bibinfo{journal}{ACM Transactions on Mathematical Software}
  \bibinfo{volume}{29}, \bibinfo{pages}{110--140}.
\bibitem[{Limkilde et~al.(2018)Limkilde, Evgrafov and
  Gravesen}]{limkilde2018reducing}
\bibinfo{author}{Limkilde, A.}, \bibinfo{author}{Evgrafov, A.},
  \bibinfo{author}{Gravesen, J.}, \bibinfo{year}{2018}.
\newblock \bibinfo{title}{On reducing computational effort in topology
  optimization: {W}e can go at least this far!}
\newblock \bibinfo{journal}{Structural and Multidisciplinary Optimization}
  \bibinfo{volume}{58}, \bibinfo{pages}{2481--2492}.
\bibitem[{Liu and Nocedal(1989)}]{liu1989limited}
\bibinfo{author}{Liu, D.C.}, \bibinfo{author}{Nocedal, J.},
  \bibinfo{year}{1989}.
\newblock \bibinfo{title}{On the limited memory {BFGS} method for large scale
  optimization}.
\newblock \bibinfo{journal}{Mathematical Programming} \bibinfo{volume}{45},
  \bibinfo{pages}{503--528}.
\bibitem[{Mavriplis(1990)}]{Mavriplis1990}
\bibinfo{author}{Mavriplis, C.}, \bibinfo{year}{1990}.
\newblock \bibinfo{title}{Proceedings of the Eighth GAMM-Conference on
  Numerical Methods in Fluid Mechanics}. \bibinfo{publisher}{Vieweg+Teubner
  Verlag}, \bibinfo{address}{Wiesbaden}. chapter \bibinfo{chapter}{A posteriori
  error estimators for adaptive spectral element techniques}.
\newblock pp. \bibinfo{pages}{333--342}.
\bibitem[{Mor{\'e} and Thuente(1994)}]{more1994line}
\bibinfo{author}{Mor{\'e}, J.J.}, \bibinfo{author}{Thuente, D.J.},
  \bibinfo{year}{1994}.
\newblock \bibinfo{title}{Line search algorithms with guaranteed sufficient
  decrease}.
\newblock \bibinfo{journal}{ACM Transactions on Mathematical Software}
  \bibinfo{volume}{20}, \bibinfo{pages}{286--307}.
\bibitem[{Nocedal and Wright(1999)}]{NW99}
\bibinfo{author}{Nocedal, J.}, \bibinfo{author}{Wright, S.J.},
  \bibinfo{year}{1999}.
\newblock \bibinfo{title}{Numerical Optimization}.
\newblock \bibinfo{publisher}{Springer-Verlag}, \bibinfo{address}{New York}.
\bibitem[{Ou and Jameson(2011)}]{ou_2011}
\bibinfo{author}{Ou, K.}, \bibinfo{author}{Jameson, A.}, \bibinfo{year}{2011}.
\newblock \bibinfo{title}{Unsteady adjoint method for the optimal control of
  advection and {B}urger's equations using high-order spectral difference
  method}, in: \bibinfo{booktitle}{49th AIAA Aerospace Sciences Meeting
  including the New Horizons Forum and Aerospace Exposition},
  \bibinfo{organization}{American Institute of Aeronautics and Astronautics}.
\bibitem[{Sagebaum et~al.(2017)Sagebaum, Albring and Gauger}]{codipack}
\bibinfo{author}{Sagebaum, M.}, \bibinfo{author}{Albring, T.},
  \bibinfo{author}{Gauger, N.R.}, \bibinfo{year}{2017}.
\newblock \bibinfo{title}{{High-Performance Derivative Computations using
  CoDiPack}}.
\newblock \bibinfo{journal}{arXiv preprint arXiv:1709.07229} .
\bibitem[{Saglietti et~al.(2016)Saglietti, Schlatter, Monokrousos and
  Henningson}]{saglietti2016adjoint}
\bibinfo{author}{Saglietti, C.}, \bibinfo{author}{Schlatter, P.},
  \bibinfo{author}{Monokrousos, A.}, \bibinfo{author}{Henningson, D.S.},
  \bibinfo{year}{2016}.
\newblock \bibinfo{title}{Adjoint optimization of natural convection problems:
  differentially heated cavity}.
\newblock \bibinfo{journal}{Theoretical and Computational Fluid Dynamics} ,
  \bibinfo{pages}{1--17}.
\bibitem[{Sandu(2006)}]{sandu2006properties}
\bibinfo{author}{Sandu, A.}, \bibinfo{year}{2006}.
\newblock \bibinfo{title}{On the properties of {R}unge-{K}utta discrete
  adjoints}, in: \bibinfo{booktitle}{International Conference on Computational
  Science}, \bibinfo{organization}{Springer}. pp. \bibinfo{pages}{550--557}.
\bibitem[{Schanen et~al.(2016)Schanen, Marin, Zhang and
  Anitescu}]{schanen_2016}
\bibinfo{author}{Schanen, M.}, \bibinfo{author}{Marin, O.},
  \bibinfo{author}{Zhang, H.}, \bibinfo{author}{Anitescu, M.},
  \bibinfo{year}{2016}.
\newblock \bibinfo{title}{Asynchronous two-level checkpointing scheme for
  large-scale adjoints in the spectral-element solver {N}ek5000}.
\newblock \bibinfo{journal}{Procedia Computer Science} \bibinfo{volume}{80},
  \bibinfo{pages}{1147--1158}.
\bibitem[{Serban and Hindmarsh(2005)}]{serban2005cvodes}
\bibinfo{author}{Serban, R.}, \bibinfo{author}{Hindmarsh, A.C.},
  \bibinfo{year}{2005}.
\newblock \bibinfo{title}{{CVODES}, the sensitivity-enabled {ODE} solver in
  {SUNDIALS}}, in: \bibinfo{booktitle}{Proceedings of the 5th International
  Conference on Multibody Systems, Nonlinear Dynamics and Control, Long Beach,
  CA}.
\bibitem[{Walther and Griewank(2009)}]{adolc}
\bibinfo{author}{Walther, A.}, \bibinfo{author}{Griewank, A.},
  \bibinfo{year}{2009}.
\newblock \bibinfo{title}{{Getting Started with ADOL-C}}, in:
  \bibinfo{booktitle}{Combinatorial {S}cientific {C}omputing}, pp.
  \bibinfo{pages}{181--202}.
\bibitem[{Zhang et~al.(2019)Zhang, Constantinescu and Smith}]{Zhang_2019}
\bibinfo{author}{Zhang, H.}, \bibinfo{author}{Constantinescu, E.M.},
  \bibinfo{author}{Smith, B.F.}, \bibinfo{year}{2019}.
\newblock \bibinfo{title}{{PETSc} {TSA}djoint: {A} discrete adjoint ode solver
  for first-order and second-order sensitivity analysis}.
\newblock \URLprefix \url{https://arxiv.org/abs/1912.07696},
  \href{http://arxiv.org/abs/1912.07696}{{\tt arXiv:1912.07696}}.
\bibitem[{Zhang and Sandu(2014)}]{zhang2014fatode}
\bibinfo{author}{Zhang, H.}, \bibinfo{author}{Sandu, A.}, \bibinfo{year}{2014}.
\newblock \bibinfo{title}{{FATODE}: a library for forward, adjoint, and tangent
  linear integration of {ODE}s}.
\newblock \bibinfo{journal}{SIAM Journal on Scientific Computing}
  \bibinfo{volume}{36}, \bibinfo{pages}{C504--C523}.

\end{thebibliography}

\newpage
{\bf Government License.}  The submitted manuscript has been created by
UChicago Argonne, LLC, Operator of Argonne National Laboratory
(``Argonne''). Argonne, a U.S. Department of Energy Office of Science
laboratory, is operated under Contract No. DE-AC02-06CH11357. The
U.S. Government retains for itself, and others acting on its behalf, a
paid-up nonexclusive, irrevocable worldwide license in said article to
reproduce, prepare derivative works, distribute copies to the public,
and perform publicly and display publicly, by or on behalf of the
Government.  The Department of Energy will provide public access to
these results of federally sponsored research in accordance with the
DOE Public Access
Plan. http://energy.gov/downloads/doe-public-access-plan.
\end{document}